\documentclass[10pt]{article}
\usepackage[utf8]{inputenc}
\usepackage[T1]{fontenc}
\usepackage[paper=a4paper,left=2.5cm,right=2.5cm,top=2.5cm,bottom=2.5cm]{geometry}
\usepackage{amsmath}
\usepackage{amssymb}
\usepackage{graphicx}
\usepackage{amsfonts}
\usepackage{graphicx}
\usepackage{epstopdf}
\usepackage{algorithmic}
\usepackage{bm}
\usepackage{subcaption}
\usepackage{amsmath}
\usepackage{amsthm}
\usepackage{cases}
\usepackage{color}
\usepackage{multirow}
\usepackage{hyperref}

\ifpdf
\DeclareGraphicsExtensions{.eps,.pdf,.png,.jpg}
\else
\DeclareGraphicsExtensions{.eps}
\fi

\newcommand{\vv}{{\bm{v}}}

\newcommand{\ff}{{\bm{f}}}
\newcommand{\qq}{{\bm{q}}}
\newcommand{\FF}{{\bm F}}

\newcommand{\HH}{{\bm H}}
\newcommand{\QQ}{{\bm Q}}
\newcommand{\LL}{{\bm L}}
\newcommand{\CC}{{\bm C}}

\newcommand{\xx}{{\bm{x}}}
\newcommand{\ppsi}{{\bm{q}}}
\renewcommand{\AA}{{\bm{A}}}
\newcommand{\BB}{{\bm B}}

\newcommand{\pdx}{{\partial_x}}
\newcommand{\pdxl}[1]{{\partial_{x_#1}}}

\newcommand{\pdt}{{\partial_t}}
\newcommand{\uu}{{\bm{u}}}
\newcommand{\nn}{{\bm{n}}}
\newcommand{\hh}{{\bm{h}}}

\newcommand{\kin}{\text{kin}}
\newcommand{\divv}{\text{div}}
\newcommand{\IMEX}{{\text{SIFV-PC} }}
\newcommand{\IM}{{\text{IMFV-PC} }}
\newcommand{\SIMHDthree}{\text{SIFV-EB} }
\newcommand{\slow}{\text{slow}}
\newcommand{\fast}{\text{fast}}
\newtheorem{theorem}{Theorem}[section]
\newtheorem{lemma}[theorem]{Lemma}

\theoremstyle{definition}

\theoremstyle{remark}
\newtheorem{remark}[theorem]{Remark}

\title{Semi-implicit relaxed finite volume schemes for hyperbolic multi-scale systems of conservation laws}
\author{Andrea Thomann$^\ast$ \\[2pt]
$^\ast$Université de Strasbourg, CNRS, Inria, IRMA, F-67000 Strasbourg, France}

\begin{document}

\maketitle

\begin{abstract}
   In this paper a new semi-implicit relaxation scheme for the simulation of multi-scale hyperbolic conservation laws based on a Jin-Xin relaxation approach is presented.
   It is based on the splitting of the flux function into two or more subsystems separating the different scales of the considered model whose stiff components are relaxed thus yielding a linear structure of the resulting relaxation model on the relaxation variables.
   This allows the construction of a linearly implicit numerical scheme, where convective processes are discretized explicitly.
   Thanks to this linearity, the discrete scheme can be reformulated in linear decoupled wave-type equations resulting in the same number of evolved variables as in the original system.
   To obtain a scale independent numerical diffusion, centred fluxes are applied on the implicitly treated terms, whereas classical upwind schemes are applied on the explicit parts.
   The numerical scheme is validated by applying it on the Toro \& V\'azquez-Cend\'on (2012) splitting of the Euler equations and the Fambri (2021) splitting of the ideal MHD equations where the flux is split in two, respectively three sub-systems.
   The performance of the numerical scheme is assessed running benchmark test-cases from the literature in one and two spatial dimensions.

\end{abstract}


\section{Introduction}

Many models that capture physical phenomena contain multiple parameters which affect the propagation of information.
In the context of hyperbolic systems, this is reflected by wave speeds that may differ several orders of magnitudes depending on the considered regimes.
In general, one can divide these dynamics into slow processes, e.g. convection or transport with the local flow velocity and fast processes, e.g. acoustic waves at large sound speeds \cite{Klein1995}, shear waves in rigid materials \cite{Dumbser2016,Brauer2016} or Alfv\'en waves in strong magnetic fields \cite{Dumbser2019,Fambri2021}.
The appearance of such different time scales, i.e. stiffness of some flux terms, is characteristic for multi-scale equations.
Their numerical resolution is posing several challenges and is thus a very active field of research.

Applying traditional explicit upwind schemes on such a multi-scale problem has several drawbacks.
Firstly, they suffer from excessive numerical diffusion in regimes where fast processes are present, thus lead to spurious solutions.
A well studied example is low Mach number flow, e.g. \cite{BoscarinoRussoScandurra2018,Dellacherie2010,GuillardViozat1999}, which can be transferred to general multi-scale systems.
Secondly, for stability, an explicit scheme requires a quite severe time step restriction which depends on the propagation speed of the fastest process.
Thus, this can lead, in extreme cases to vanishing time steps, but at least to an elevated number of time steps compared to what slow dynamics would require.
Being forced to use such small time steps results into unnecessary computational overhead.
As a side effect all waves included in the model are resolved even though at times only slow convective processes are of interest and inaccuracies and spurious artefacts on the fast waves are acceptable.

Therefore, implicit or partially implicit solvers became of increased interest, see e.g. \cite{Oeffner2025,AbbIolPup2017,Bispen2017,Boscarino2016,CoulFraHelRatSon2019} and references therein.
Their use allows fast processes to be integrated implicitly which means the time step can be oriented towards slower processes without the loss of stability.
Moreover, using a central discretization for implicitly treated stiff terms mitigates the excessive numerical diffusion which would appear when an upwind discretization were used for the discretization of fast processes.

However, due to the non-linearity of the underlying models, non-linear, often coupled, systems of equations need to be solved implicitly which require expensive iterative solvers.
Especially in the case of multi-scale problems with large wave speeds the resulting implicit systems are generally ill conditioned thus the convergence cannot be guaranteed or is very slow.

Thus a discretization which results in linear implicit systems are preferred as they can be solved with direct solvers or fast iterative methods equipped with standard preconditioners to take care of ill conditioned systems.
Such a linearization can be achieved by either modifying and linearizing the equations directly\cite{Bispen2017}, making use of linearities inside the model in certain configurations \cite{BoscarinoRussoScandurra2018,Degond2011}, or utilizing semi-implicit techniques \cite{Boscarino2016, Boscheri2024} where a suitable linearization in time is applied on the non-linear terms.
However, some of these approaches are quite intrusive and model specific.
Thus the derivation of the numerical scheme needs to be investigated anew for each problem under consideration.

General approaches are quite rare. In \cite{AbbIolPup2017,Thomann2023} unconditionally stable linearly implicit methods for general hyperbolic systems of conservation laws were derived.
In particular, in \cite{Thomann2023} a Jin-Xin relaxation model \cite{JinXin1995} was the basis for the development of the numerical scheme and the resulting discretized equations were rewritten into linear implicit wave-type equations.
This was advantageous in two aspects, avoiding the evolution of additional relaxation variables and decoupled linear systems with positive definite coefficient matrices with a Laplacian structure.
The scheme showed improved properties on the resolution of slow moving contact waves compared to linearized models due to a prediction correction procedure projecting the solution obtained from frozen characteristics on the original flux dynamics.

However, on long time scales this so-called \textit{relaxed fully implicit scheme} \cite{Thomann2023} proved to be quite diffusive which was in particular observed in transport processes where fast acoustic waves were present.
Thus a reduction of the time step is necessary to capture the slow scales with sufficient accuracy.
Motivated by these findings, in the current work, a splitting of the flux which drives the evolution of the variables is performed before a Jin-Xin relaxation is applied on the non-linear stiff flux part which triggers the fast scales.
This allows to treat the slow convective or transport processes explicitly, while fast waves are integrated implicitly.
Entering the framework of \cite{Thomann2023}, the resulting implicit systems are formulated as decoupled wave-type equations yielding similar computational effort as the fully implicit scheme \cite{Thomann2023}.
Moreover, in the case of more than one fast scale, several implicit subsystems can be considered.
This is for instance the case for the ideal Magneto-Hydro-Dynamic (MHD) equations where besides the Mach number which scales the sound waves, the Alfv\'en number appears which is induced by the strength of the magnetic field.
Therefore, splitting the flux in three sub-fluxes as done in \cite{Fambri2021,Boscheri2024} allows a separation of the scales.
Together with the splitting of the Euler equations in two sub-fluxes introduced by Toro \& V\'azquez-Cend\'on \cite{Toro2012}, the splitting from Fambri \cite{Fambri2021} for the MHD equations will be considered to test and illustrate the new semi-implicit finite volume scheme.

The paper is organised as follows.
In Section \ref{sec:Problem}, the problem is described and the concept of Jin-Xin relaxation is recalled and applied on the flux splitting.
Subsequently, in Section \ref{sec:Scheme}, the numerical scheme is derived based on a finite volume framework.
First, the numerical method is illustrated deriving a first order semi-implicit scheme based on the splitting in one explicitly and one implicitly treated sub-flux which is then extended to two implicitly treated sub-fluxes.
Then the high order extension for the two-split semi-implicit scheme is discussed and three-split schemes using partitioned Runge-Kutta methods.
Finally the multi-dimensional numerical scheme is given, before the schemes are illustrated at the examples of the Euler and ideal MHD equations.
In Section \ref{sec:NumRes}, the scheme is applied on a series of benchmark tests from the literature for the Euler and ideal MHD equations.
This includes the assessment of accuracy, Riemann Problems and two-dimensional flows in varying Mach and Alfv\'en number regimes.
Moreover, the asymptotic preserving and contact property of the scheme for the Euler equations proven in Section \ref{sec:Scheme} are numerically verified.
A section of conclusion and outlook to future developments concludes the paper.

\section{Problem description and relaxation models}
\label{sec:Problem}
For the sake of clarity in the derivations, in this section, we consider a system of hyperbolic conservation laws in one space dimension (1D) given by
\begin{equation}\label{sys:CL}
    \pdt \qq + \pdx \ff(\qq) = 0,
\end{equation}
where $\qq \in \mathbb{R}^m$ is the state vector of conserved variables and $\ff: \mathbb{R}^m \to \mathbb{R}^m$ a flux function.
Since we consider a hyperbolic system, we can decompose the Jacobian $\nabla_\qq \ff(\qq) = R \Lambda R^{-1}$ where the characteristic speeds $\lambda_i(\qq), i = 1, \dots, m$ form the diagonal matrix $\Lambda$.

In \cite{JinXin1995}, the authors proposed a linearization of \eqref{sys:CL} in terms of a relaxation technique, later on called Jin-Xin relaxation, where the whole flux $\ff$ was linearized as follows
\begin{subequations}\label{sys:Relax}
    \begin{align}
        \pdt \ppsi + \pdx \vv &=0 \\
        \pdt \vv + \AA^2 \pdx \ppsi &= \frac{1}{\varepsilon}(\ff(\ppsi)- \vv).
    \end{align}
\end{subequations}
Therein $\varepsilon > 0$ denotes the relaxation rate and in the diagonal matrix $\AA$ a constant approximation of the characteristic speeds is stored.
In \cite{JinXin1995}, it was shown under sufficient smoothness conditions and for small $\varepsilon > 0$ that the original system \eqref{sys:CL} can be recovered in the limit $\varepsilon \to 0$ under the sub-characteristic condition
\begin{equation}
    \AA^2 - \left(\nabla_\qq \ff(\qq)\right)^2 \geq 0 \quad \text{(positive semi-definite)}.
\end{equation}
In practise setting
\begin{equation}\label{eq:Aset1D}
    \AA = a \mathbb{I} \in \mathbb{R}^{m\times m}
\end{equation}
where the value $a \in \mathbb{R}$ is the constant approximation $a = \max_{i=1,\dots,m} |\lambda_i(\qq)|$ satisfies the sub-characteristic condition and requires only the knowledge of the maximal characteristic speed.
Thus the knowledge of the maximal absolute eigenvalue suffices to ensure the formal limit to the original system \eqref{sys:CL}.

For systems, where all characteristic speeds are close together, to approximate all wave speeds with the maximal one, does only introduce minimal diffusion.
However, in systems where the characteristic speeds differ several order of magnitudes, the slow speeds are not well approximated.
This structure is reflected in dynamics where $\ff$ can be decomposed into two parts $\ff = \ff^{\slow} + \ff^{\fast}$ where $\ff^\slow$ describes slow dynamics, i.e. is associated to slow characteristic speeds, and $\ff^\fast$ with characteristic speeds which are several orders faster than the slow ones.
Moreover, we require that $\ff^{\fast}$ has real eigenvalues and is diagonalisable and $\ff^{\slow}$ has real eigenvalues where the associated sub-system is (at least) weakly hyperbolic.
Then the system reads
\begin{equation}\label{sys:CLsplit}
    \pdt \ppsi + \pdx \ff^{\slow}(\ppsi) + \pdx \ff^{\fast}(\ppsi) = 0,
\end{equation}
Since $\ff^{\slow}$ is non-stiff and will be treated in the numerical scheme later on explicitly, its linearization through relaxation is not necessary.
However, $\ff^\fast$ contains stiff terms and will be treated implicitly in the numerical scheme.
Therefore, to avoid non-linear implicit systems, we relax only the non-linear flux $\ff^{\fast}$ using the Jin-Xin relaxation approach as described above.
We approximate $\ff^\fast$ with a new relaxation variable $\vv$ and obtain
\begin{subequations}\label{sys:RelaxSplit}
    \begin{align}
        \pdt \ppsi + \pdx \ff^{\slow}(\ppsi) + \pdx \vv &= 0 \\
        \pdt \vv + \AA_\fast^2 \pdx \ppsi &= \frac{1}{\varepsilon}(\ff^{\fast}(\ppsi)- \vv).
    \end{align}
\end{subequations}
Therein $\AA_\fast$ is now a diagonal matrix which contains a constant approximation of the characteristic speeds mainly associated with $\ff^\fast$.
Applying a Chapman-Enskog expansion for small $\varepsilon > 0$ and sufficiently smooth $\ff^{\fast}$, we have the following result.
\begin{lemma}
    The relaxation system \eqref{sys:RelaxSplit} is dissipative under the sub-characteristic condition $\AA_\fast^2 - (\nabla_\ppsi \ff^\text{fast} )^T (\nabla_\ppsi\ff) > 0$ and $\varepsilon > 0$ and formally tends to the original system \eqref{sys:CL} as $\varepsilon \to 0$.
\end{lemma}
Analogously to \eqref{eq:Aset1D}, we set in practise
\begin{equation}
    \AA = a^\fast ~\mathbb{I} \in \mathbb{R}^{m\times m}
\end{equation}
where the value $a^\fast$ is the constant approximation $a^\fast > \max_i |\lambda_i^\fast(\qq)|$, where $\lambda_i^\fast$ are the eigenvalues of $\nabla_\qq \ff^\fast(\qq)$ such that the sub-characteristic condition is fulfilled.

Indeed the result can be extended to a decomposition into more than two subsystems.
Let $\ff = \ff^\slow + \sum^k_{j=1} \ff^\fast_j$ where the flux $\ff$ is split into $\ff^\slow$ describing the slow dynamics and $k$ fluxes $\ff^\fast_j$ associated to fast processes whose stiffness is influenced by a different model parameter, e.g. the Mach number, Alfv\'en number etc. respectively.
As before, we assume that each flux $\ff_j^\fast$ is diagonalisable with real eigenvalues and $\ff^\slow$ has real eigenvalues and the associated sub-system is (at least) weakly hyperbolic.
Introducing now $k$ new variables $\vv_j, j = 1, \dots, k$ which approximate the fast fluxes $\ff^\fast_j$ respectively, we can write the following relaxation system using Jin-Xin relaxation
\begin{subequations}\label{sys:RelaxMultiSplit}
    \begin{align}
        \pdt \ppsi + \pdx \ff^{ex}(\ppsi) + \sum_{j=1}^{k}\pdx \vv_j &= 0 \\
        \pdt \vv_j + \AA_{j}^2 \pdx \ppsi &= - \frac{1}{\varepsilon}(\ff_j^{im}(\ppsi)- \vv_j), \quad j=1, \dots, k
    \end{align}
\end{subequations}
Applying Chapman-Enskog expansion, we obtain the following result
\begin{lemma}\label{lem:subchar_multisplit_1D}
    The relaxation system \eqref{sys:RelaxMultiSplit} is dissipative under the sub-characteristic conditions $\AA_{j}^2 - (\nabla \ff_j^\fast)^T (\nabla_\qq \ff) > 0$ for all $j = 1, \dots, k$ and $\varepsilon >0$ and formally tends to the original system \eqref{sys:CL} as $\epsilon \to 0$.
\end{lemma}
In practise, we set
\begin{equation}
    \AA_j = a_j^\fast ~\mathbb{I} \in \mathbb{R}^{m\times m}
\end{equation}
where the value $a_j^\fast$ is the constant approximation $a_j^\fast > \max_i |\lambda_i^\fast(\qq)|$, where $\lambda_i^\fast$ are the eigenvalues of $\nabla_\qq \ff_j^\fast(\qq)$ such that the sub-characteristic conditions are fulfilled.

An extension to two or higher space dimensions is straightforward applying the relaxation technique dimension by dimension as done in the original work of Jin\& Xin \cite{JinXin1995}.
In the following, we shortly recall the two-dimensional (2D) Jin-Xin relaxation.
Consider a hyperbolic system of conservation laws in 2D with $\xx = (x_1, x_2)^T$ given by
\begin{equation}\label{sys:CL2D}
    \pdt \ppsi + \pdxl{1} \ff^1(\ppsi) + \pdxl{2} \ff^{2}(\ppsi)= 0.
\end{equation}
Therein the fluxes $\ff^1$ and $\ff^2$ denote the fluxes in $x_1$- and $x_2$- direction respectively.
Introducing two new sets of variables $\vv^1, \vv^2 \in \mathbb{R}^m$, we can write the following relaxation system
\begin{subequations}\label{sys:Relax2D}
    \begin{align}
        \pdt \ppsi + \pdxl{1} \vv^1 + \pdxl{2}\vv^2 &=0 \\
        \pdt \vv^1 + \AA_1^2 \pdxl{1} \ppsi &= \frac{1}{\varepsilon}(\ff^1(\ppsi)- \vv^1)\\
        \pdt \vv^2 + \AA_2^2 \pdxl{2} \ppsi & = \frac{1}{\varepsilon}(\ff^2(\ppsi) - \vv^2),
    \end{align}
\end{subequations}
where $\AA_l, l=1,2$ are constant diagonal matrices.
Then in \cite{JinXin1995} it was shown that the relaxation system is dissipative for $\varepsilon > 0$ under the following sub-characteristic condition
\begin{equation}
    \delta_{lr} \AA_l^2 - \nabla_\qq \ff^l \nabla_\qq \ff^r > 0, \quad l,r = 1,2.
\end{equation}
Therein, $\delta_{lr}$ denotes the Kronecker delta.

Now we apply a flux splitting in the analogous way as done in the one dimensional case as follows
\begin{equation}
    \ff^{l} = \ff^{l,\slow} + \sum^k_{j=1} \ff^{l,\fast}_j, \quad l = 1,2.
\end{equation}
Introducing $2k$ sets of new variables $\vv_j^l \in \mathbb{R}^m, j=1, \dots, k,~ l = 1, 2$, we can write the following relaxation system
\begin{subequations}\label{sys:RelaxMultiSplit_2D}
    \begin{align}
        \pdt \ppsi + \sum_{l=1}^2\pdxl{l} \ff^{l,\slow}(\ppsi) + \sum_{l=1}^2\sum_{j=1}^{k}\pdxl{l} \vv_j^l &= 0 \\
        \pdt \vv_j^l + \AA_{l,j}^2 \ \pdxl{l} \ppsi &= \frac{1}{\varepsilon}(\ff_j^{l,\fast}(\ppsi)- \vv_j^l), \quad j=1, \dots, k, \quad l = 1,2.
    \end{align}
\end{subequations}
Analogously to the one dimensional case, the matrices $\AA_{l,j}$ are constant diagonal matrices and the relaxation system \eqref{sys:RelaxMultiSplit_2D} is dissipative under the sub-characteristic condition
\begin{equation}\label{eq:subchar_multisplit_2D}
    \delta_{lr} \AA_{l,j}^2 - (\nabla_\qq \ff^{l,\fast}_j)^T (\nabla_\qq \ff^r) > 0, \quad l,r = 1,2, \quad j=1,\dots,k.
\end{equation}
This framework allows a straightforward extension also to three spatial dimensions which is thus omitted here.
We conclude this section with some remarks.
\begin{remark}
    \verb| |
    \begin{itemize}
    \item Note that the dimensionality of the relaxation system increases with each added spatial dimension and fast sub-system $\ff_j^\fast$.
    Therefore, the aim of the construction of the numerical scheme is to reduce as much as possible the number of equations that need to be solved as well as the number of variables that need to be updated.
    In particular this means the construction of \textit{relaxed} schemes which are independent of $\varepsilon$ and the relaxation variables $\vv_j^l$.
    \item To reduce the dissipation introduced by the relaxation process, in practise we want to choose the smallest diagonal matrix $\AA_{l,j}$ that satisfies the respective sub-characteristic condition oriented closely towards the characteristic speeds of the original system or sub-systems. The choice that is used to perform the numerical test cases is given for each case in Section \ref{sec:NumRes}.
    \end{itemize}
\end{remark}
\section{Numerical scheme}
\label{sec:Scheme}
To construct a numerical scheme, we rely on a finite volume discretization in space, where the state variables $\ppsi$ are discretized on a primal grid and the approximated fluxes represented by the relaxation variables $\vv$ are discretized on a staggered grid since the fluxes are evaluated at the interfaces of the primal grid cells.
For clarity in derivation, we focus on the 1D case. An extension to 2D is given at the end of Section \ref{sec:Scheme}.

Let the computational domain in space be given by an interval $I\subset\mathbb{R}$ which is divided in $N$ cells $C_i$, $i=1, \dots, N$ which are defined as $C_i = (x_{i-1/}, x_{i+1/2})$ with cell center $x_i = \frac{x_{i-1/2} + x_{i+1/2}}{2}$ and $x_{i+1/2} = N \Delta x + \frac{\Delta x}{2}$ with uniform gridsize $\Delta x$.
This constitutes the primal grid.
The cell averages of the state variables $\qq$ on $C_i$ at a time $t$ are approximated by
\begin{equation}
    \qq_i(t) \approx \frac{1}{\Delta x} \int_{C_i} \qq(x,t) dx.
\end{equation}
Further we consider a dual grid shifted by $\frac{\Delta x}{2}$, i.e. the cells are given by $D_i = (x_i, x_{i+1})$ with cell center $x_{i+1/2}$ on which the relaxation variables $\vv$ will be discretized.

The time-line is discretized by $t^n = n \Delta t$, where $\Delta t$ is a time increment yet to be determined for the numerical scheme to ensure stability.

We begin with overall first order schemes.

\subsection{First order semi-implicit schemes}
For clarity in the construction steps of the numerical scheme, we consider first system \eqref{sys:CLsplit} with a splitting of the flux into two associated sub-systems.
The associated relaxation system which forms the basis of the numerical scheme is given by system \eqref{sys:RelaxSplit}.
Note that the evolution equation for the state variable $\qq$ is linear in $\vv$ and the vice versa the equation for the relaxation variable $\vv$ is linear in $\qq$.
Moreover, the non-linear flux $\ff^\fast$ appears in the source term which is free of any derivatives, i.e. it is algebraic.
To make use of this structure, we apply an operator splitting separating the algebraic source term from the homogeneous left hand side.
Since we are interested in a relaxed numerical scheme, i.e. the limit $\varepsilon = 0$, we will make use of the exact relation $\vv = \ff^\fast(\ppsi)$, i.e. project the solution on the equilibrium manifold for $\varepsilon = 0$ which serves as corrector later on in the numerical scheme.

The homogeneous system
\begin{subequations}\label{sys:RelaxSplit_hom}
    \begin{align}
        \pdt \ppsi + \pdx \ff^{\slow}(\ppsi) + \pdx \vv &= 0, \\
        \pdt \vv + \AA_\fast^2 \pdx \ppsi &= 0,
    \end{align}
\end{subequations}
will yield the predicted solution in the numerical scheme denoted by $\qq^{(1)}, \vv^{(1)}$.
To avoid severe time-step restrictions due to the (potentially) stiff sub-flux $\ff^\fast$ which is approximated by $\vv$, the relaxation variables $\vv$ will be integrated implicitly while the sub-flux associated to the slow dynamics $\ff^\slow$ is treated in an explicit manner.
A first order discretization is realized by evaluating the state variables $\qq_i$ on the primal and the relaxation variables $ \vv_{i+1/2}$ on the dual grid and using an implicit-explicit Euler time integration method.
Since in \eqref{sys:RelaxSplit_hom} the characteristics associated to the fast dynamics are "frozen" in time due to the linearization, a numerical solution $\qq_i^{(1)}, \vv_{i+1/2}^{(1)}$ of \eqref{sys:RelaxSplit_hom} needs to be corrected projecting on the equilibrium manifold, i.e. utilizing the relation $\vv = \ff^\fast(\qq)$.
For more details and illustration on the effect of the final corrected update, see \cite{Thomann2023}.

The numerical scheme with initialization $\vv^n = \ff^\fast(\qq^n)$ is then given by
\begin{subequations}\label{sys:RelaxSplit_disc}
    \begin{eqnarray}
        \ppsi^{(1)}_i - \ppsi^n_i+ \frac{\Delta t}{\Delta x}(\FF^{ex,n}_{i+1/2} - \FF^{ex,n}_{i-1/2}) + \frac{\Delta t}{\Delta x}(\vv_{i+1/2}^{(1)} - \vv_{i-1/2}^{(1)}) &= &0, \label{eq:Relax_Split_disc_psi} \\
        \vv_{i+1/2}^{(1)} - \vv_{i+1/2}^{n} + \frac{\Delta t}{\Delta x}\AA_\fast^2 ( \ppsi_{i+1}^{(1)} - \ppsi_i^{(1)}) &= & 0,\label{eq:Relax_Split_disc_v} \\
        \ppsi^{n+1} = \ppsi^n - \frac{\Delta t}{\Delta x}\left((\FF^{ex,n}_{i+1/2} - \FF^{ex,n}_{i-1/2})+ (\FF^{im,(1)}_{i+1/2} - \FF^{im,(1)}_{i-1/2})\right).&& \label{eq:Relax_Split_disc_update}
    \end{eqnarray}
\end{subequations}
The numerical flux $\FF^{ex,n}_{i+1/2}$ is a Rusanov flux based on the characteristic speeds of the slow dynamics $\ff^\slow$
\begin{equation}\label{eq:flux_ex_Rusanov}
    \FF_{i+1/2}^{ex}(\ppsi_i, \ppsi_{i+1}) = \frac{1}{2} (\ff^{ex}_i + \ff^{ex}_{i+1}) - \frac{1}{2} \max_j(|\lambda^\slow_j(\ppsi_i)|,|\lambda^\slow_{j}(\ppsi_{i+1})|)(\ppsi_{i+1} - \ppsi_i)
\end{equation}
where $\lambda_j^\slow$ denote the characteristic speeds associated to $\ff^\slow$.
The approximated fluxes $\vv$ are discretized with a central scheme
\begin{equation}
    \vv_{i+1/2} = \frac{1}{2}\left(\vv_i + \vv_{i+1}\right),
\end{equation}
and the implicit numerical fluxes are discretized with central differences
\begin{equation}\label{eq:flux_im_center2}
    \FF_{i+1/2}^{im} = \frac{1}{2} (\ff^{\fast}(\qq_i) + \ff^\fast(\qq_{i+1})).
\end{equation}
Thus, a time step restriction of Rusanov-type for stability has to be imposed
\begin{equation}\label{eq:CFL}
    \Delta t \leq \nu \frac{\Delta x}{\max_j|\lambda^\slow_j|}, \quad 0 < \nu \leq 1.
\end{equation}

Due to the relaxation procedure, the number of variables has doubled.
To reduce the resulting computational overhead of evaluating the coupled implicit, but linear, systems \eqref{eq:Relax_Split_disc_psi},\eqref{eq:Relax_Split_disc_v}, we replace the implicit fluxes $\vv_{i+1/2}^{(1)}$ in \eqref{eq:Relax_Split_disc_psi} by their evolution given in \eqref{eq:Relax_Split_disc_v}.
This yields the following wave equation for the state variables
\begin{equation}\label{eq:impl_sys_imex}
    \ppsi^{(1)}_i - \frac{\Delta t^2}{\Delta x^2} \AA_\fast^2 \left(\ppsi_{i+1}^{(1)} - 2 \ppsi_i^{(1)} + \ppsi_{i-1}^{(1)}\right)= \ppsi^n_i - \frac{\Delta t}{\Delta x}(\FF^{ex,n}_{i+1/2} - \FF^{ex,n}_{i-1/2}) - \frac{\Delta t}{\Delta x} (\FF^{im,n}_{i+1/2} - \FF^{im,n}_{i-1/2}).
\end{equation}
Note that due to $\AA_\fast$ being a constant diagonal matrix, the implicit system \eqref{eq:impl_sys_imex} decouples into $m$ linear systems of dimension $N$, where $m$ is the number of state variables and $N$ the number of discretization cells.
Moreover, the coefficient matrix is symmetric and positive definite which allows the use of efficient linear solvers.
Summarizing, the first order two-split numerical scheme is given by the linear implicit systems \eqref{eq:impl_sys_imex} followed by the final update \eqref{eq:Relax_Split_disc_update}.

Next, we consider an extension of the two-split scheme \eqref{eq:impl_sys_imex}\eqref{eq:Relax_Split_disc_update} to a three-split scheme, i.e. for the relaxation system \eqref{sys:RelaxMultiSplit} with $k=2$ fast sub-fluxes $\ff_1^\fast, \ff_2^\fast$.
A splitting in more than one fast part, is useful if characteristic speeds with different orders of magnitude appear in the model under consideration.

Analogously to the derivation of the two-split scheme, the flux $\ff^\slow$ is treated explicitly, while the remaining terms are treated implicitly.
Thus, a first-order prediction-correction numerical scheme is given by
\begin{subequations}\label{sys:RelaxSplit3_disc}
    \begin{eqnarray}
        &&\ppsi^{(1,\ast)}_i - \ppsi^n_i+ \frac{\Delta t}{\Delta x}(\FF^{ex,n}_{i+1/2} - \FF^{ex,n}_{i-1/2}) + \frac{\Delta t}{\Delta x}(\vv_{1,i+1/2}^{(1,\ast)} - \vv_{1,i-1/2}^{(1,\ast)}) = 0, \label{eq:Relax_Split3_disc_psi} \\
        &&\vv_{1,i+1/2}^{(1,\ast)} - \vv_{1,i+1/2}^{n} + \frac{\Delta t}{\Delta x}\AA_{1}^2 ( \ppsi_{i+1}^{(1,\ast)} - \ppsi_i^{(1,\ast)}) =  0,\label{eq:Relax_Split3_disc_v1} \\
        &&\ppsi^{(1)}_i - \ppsi^{(1,\ast)}_i+ \frac{\Delta t}{\Delta x}(\vv_{2,i+1/2}^{(1,\ast)} - \vv_{2,i-1/2}^{(1,\ast)}) = 0, \label{eq:Relax_Split3_disc_psi1} \\
        &&\vv_{2,i+1/2}^{(1)} - \vv_{2,i+1/2}^{(1,\ast)} + \frac{\Delta t}{\Delta x}\AA_{2}^2 ( \ppsi_{i+1}^{(1)} - \ppsi_i^{(1)}) =  0,\label{eq:Relax_Split3_disc_v2} \\
        &&\ppsi^{n+1} = \ppsi^n - \frac{\Delta t}{\Delta x}\left((\FF^{ex,n}_{i+1/2} - \FF^{ex,n}_{i-1/2})+ (\FF^{im,(1)}_{1,i+1/2} - \FF^{im,(1)}_{1,i-1/2}) + (\FF^{im,(1)}_{2,i+1/2} - \FF^{im,(1)}_{2,i-1/2})\right). \label{eq:Relax_Split3_disc_update}
    \end{eqnarray}
\end{subequations}
Therein, w.l.o.g. first the dynamics with respect to $\ff^\fast_1$ expressed by the relaxation variables $\vv_{1}$ are integrated implicitly yielding the intermediate updates $\qq_i^{(1,\ast)}, \vv^{(1,\ast)}_{1,i+1/2}$ from equations \eqref{eq:Relax_Split3_disc_psi}, \eqref{eq:Relax_Split3_disc_v1}.
Then, using this information, the dynamics with respect to $\ff^\fast_2$ expressed by the relaxation variables $\vv_2$ are integrated implicitly yielding the intermediate states $\qq_i^{(1)}, \vv^{(1)}_{2,i+1/2}$ which concludes the prediction step.
Then, using the final predicted variables, the correction is given in \eqref{eq:Relax_Split3_disc_update} making use of the projection $\vv_1 = \ff_1^\fast(\qq)$ and $\vv_2 = \ff^\fast_2(\qq)$.
The numerical flux $\FF_{i+1/2}^{ex}$ is defined as in \eqref{eq:flux_ex_Rusanov}, and the fluxes $\FF_j^{im,\bullet}, j=1,2$ are centred fluxes defined by
\begin{equation}
    \FF_{j,i+1/2}^{im,\bullet} = \frac{1}{2} (\ff_j^{\fast}(\qq_i^\bullet) + \ff_j^\fast(\qq_{i+1}^\bullet)), \quad j = 1,2.
\end{equation}
Therefore, the CFL condition of the two-split scheme \eqref{eq:CFL} holds also for the three-split scheme.

Analogously to the two-split scheme, we reduce the number of updated variables by replacing $\vv_{j,i+1/2}$ in \eqref{eq:Relax_Split3_disc_psi} and \eqref{eq:Relax_Split3_disc_psi1} by their evolution given in \eqref{eq:Relax_Split3_disc_v1} and \eqref{eq:Relax_Split3_disc_v2} respectively which yields the following implicit systems
\begin{subequations}\label{eq:impl_sys_imex_3}
    \begin{align}
    \ppsi^{(1,\ast)}_i - \frac{\Delta t^2}{\Delta x^2} \AA_{1}^2 \left(\ppsi_{i+1}^{(1,\ast)} - 2 \ppsi_i^{(1,\ast)} + \ppsi_{i-1}^{(1,\ast)}\right) &= \ppsi^n_i - \frac{\Delta t}{\Delta x}(\FF^{ex,n}_{i+1/2} - \FF^{ex,n}_{i-1/2}+\FF^{im,n}_{1,i+1/2} - \FF^{im,n}_{1,i-1/2})\label{eq:impl_sys_imex_3_s1star}\\
    \ppsi^{(1)}_i - \frac{\Delta t^2}{\Delta x^2} \AA_{2}^2 \left(\ppsi_{i+1}^{(1)} - 2 \ppsi_i^{(1)} + \ppsi_{i-1}^{(1)}\right)& = \ppsi_i^{(1,\ast)} - \frac{\Delta t}{\Delta x} (\FF^{im,(1,\ast)}_{1,i+1/2} - \FF^{im,(1,\ast)}_{1,i-1/2}). \label{eq:impl_sys_imex_3_s1}
    \end{align}
\end{subequations}
Therein we have used that at time $t^n$, the data is in relaxation equilibrium, i.e. it holds $\vv_1^n = \ff_1^{im}(\ppsi^n)$ and $\vv_2^n = \ff_2^{im}(\ppsi^n)$.
Summarizing, the three split first order semi-discrete scheme is given by the implicit systems \eqref{eq:impl_sys_imex_3} and the update \eqref{eq:Relax_Split3_disc_update}.
As remarked for the two-split scheme, the implicit systems decouple into $2m$ systems of dimension $N$. Since in general $\AA_1 \neq \AA_2$, the coefficient matrices differ for the two intermediate stages, however they are decoupled and have a favourable sparse banded structure.

Along these lines, a framework for the construction of first order numerical schemes with even more than two implicit fluxes can be constructed.
The resulting scheme will be first order in time, see \cite{Hairer1991} for a discussion on splitting methods.
However, the quality of the numerical results and applicability of the presented numerical schemes depend on a suitable and robust flux splitting in the considered hyperbolic model taking into account the physical dependencies and couplings between the flux terms.

\subsection{High order two-split schemes}
Staying in the one-dimensional setting, we can extend the first order semi-implicit two-split scheme \eqref{eq:impl_sys_imex}, \eqref{eq:Relax_Split_disc_update} to high order of accuracy using the IMEX Runge-Kutta framework combined with a high order reconstruction procedure for the interface values in space based on the method of lines.
Here, we will use exemplarily as proof of concept a second order space discretization which can easily be extended using high order discretization of the Laplacian based on finite differences and high order reconstruction at the cell interfaces, e.g. WENO or ENO reconstructions \cite{Cravero2018,Shu1998}.

We start again from the relaxation model \eqref{sys:RelaxSplit} with an explicitly treated flux $\ff^\slow$ and implicit integration of the fast dynamics associated to $\vv$. Using a diagonally implicit IMEX-RK scheme with $s$ stages and Butcher tableaux given for the explicit part by $((\tilde \alpha)_{r,k=1}^{s}, (\tilde \beta)_{r=1}^s)$ and the implicit part by $(\alpha_{r,k=1}^{s}, \beta_{r=1}^s)$, we obtain the following approximation for $r=1, \dots, s$
\begin{subequations}\label{sys:RelaxSplit_disc_ho}
    \begin{eqnarray}
        \ppsi^{(r)}_i - \ppsi^n_i+ \sum_{k=1}^{r-1}\tilde{\alpha}_{rk}\frac{\Delta t}{\Delta x}(\FF^{ex,n}_{i+1/2} - \FF^{ex,n}_{i-1/2}) + \sum_{k=1}^{r}\alpha_{rk} \frac{\Delta t}{\Delta x}(\vv_{i+1/2}^{(k)} - \vv_{i-1/2}^{(k)}) &= &0, \label{eq:Relax_Split_disc_psi_ho} \\
        \vv_{i+1/2}^{(r)} - \vv_{i+1/2}^{n} + \sum_{k=1}^{r}\alpha_{rk}\frac{\Delta t}{\Delta x}\AA_\fast^2 ( \ppsi_{i+1}^{(k)} - \ppsi_i^{(k)}) &= & 0,\label{eq:Relax_Split_disc_v_ho} \\
        \ppsi^{n+1} = \ppsi^n - \frac{\Delta t}{\Delta x}\sum_{r=1}^{s} \left(\tilde \beta_{r}(\FF^{ex,(r)}_{i+1/2} - \FF^{ex,(r)}_{i-1/2})+ \beta_r(\FF^{im,(r)}_{i+1/2} - \FF^{im,(r)}_{i-1/2})\right).&& \label{eq:Relax_Split_disc_update_ho}
    \end{eqnarray}
\end{subequations}
Therein, the numerical fluxes are given by \eqref{eq:flux_ex_Rusanov}, \eqref{eq:flux_im_center2}.
Then, with the same procedure as in the first order case, we replace $\vv^{(r)}$ in \eqref{eq:Relax_Split_disc_psi_ho} by the stage discretization \eqref{eq:Relax_Split_disc_v_ho}.
At the same time we can make use of the projection, i.e. that computing the $r$-th stage, the projection yields $v^{(k)} = \ff^{(im)}(\ppsi^{(k)})$ for $k = 1, \dots, r-1$.
Since we start from data in relaxation equilibrium, we have in addition $\vv^n = \ff^{im}(\ppsi^n)$.
Then we obtain for the $r$-th stage
\begin{align}\label{eq:impicit_2split_linsys_ho}
    \begin{split}
        \ppsi^{(r)}_i &- \frac{\alpha_{rr}^2\Delta t^2}{\Delta x^2}A_{im}^2 \left(\ppsi_{i+1}^{(r)} - 2 \ppsi_i^{(r)} + \ppsi_{i-1}^{(r)}\right) = \\
        \ppsi^n_i &- \sum_{k=1}^{r-1}\left(\tilde{\alpha}_{rk}\frac{\Delta t}{\Delta x}(\FF^{ex,n}_{i+1/2} - \FF^{ex,n}_{i-1/2}) + \sum_{k=1}^{r}\alpha_{rk} \frac{\Delta t}{\Delta x}(\FF_{i+1/2}^{im,(k)} - \FF_{i-1/2}^{im,(k)})\right) \\
        &+\alpha_{rr}\left(\FF_{i+1/2}^{im,n}- \FF_{i-1/2}^{im,n} + \frac{\Delta t}{\Delta x}\AA_{im}^2\sum_{k=1}^{r-1}\alpha_{rk} ( \ppsi_{i+1}^{(k)} - 2\ppsi_i^{(k)} + \ppsi_{i-1}^{(k)})\right)
    \end{split}
\end{align}
For a single diagonal implicit RK method (SDIRK), the $\alpha_{rr}$ coefficients coincide and the coefficient matrix in \eqref{eq:impicit_2split_linsys_ho} is the same for all stages and variables.
Moreover, since $\AA_{\fast}$ is a diagonal matrix with constant entries the system decouples into $m$ sub-systems for $\ppsi \in \mathbb{R}^m$ thus in total $s~m$ implicit systems of dimension $N$ have to be solved.

In case of the second order space discretization, the flux disretization $\FF_{i+1/2}^{im}$ is already second order accurate.
In the explicit numerical flux $\FF_{i+1/2}^{ex}$ a reconstruction of the state variables at the interface is necessary.
In case of smooth solutions, we have apply a linear reconstruction, whereas for Riemann problems and non-smooth test cases we apply a minmod limiter to obtain the interface values $\ppsi_{i+1/2}^\pm$ which are then the arguments of the Rusanov flux \eqref{eq:flux_ex_Rusanov}.
For details on the reconstruction we refer the reader to the e.g.  \cite{Toro2009}.
\subsection{Three-split schemes based on partitioned Runge-Kutta methods}
In case of the three-split scheme, the formalism of IMEX-RK methods cannot be applied.
Therefore we use the so-called partitioned or semi-implicit Runge-Kutta methods \cite{Hairer1991,Boscarino2016} which recently have been applied in the context of a three-split scheme for the MHD equations \cite{Boscheri2024}.
A three-split or multi-split relaxation system \eqref{sys:RelaxMultiSplit} can be written in the formalism of partitioned Runge-Kutta methods \cite{Boscarino2016}
\begin{subequations}\label{eq:semi-implicit}
\begin{align}
    \pdt \QQ^{ex} = -\pdx\HH(\QQ^{ex}(t), \QQ^{im}(t)) = -\pdx\HH_1(\QQ^{ex}(t)) - \pdx\HH_2(\QQ^{im}(t)), \\
    \pdt \QQ^{im} = -\pdx\HH(\QQ^{ex}(t), \QQ^{im}(t)) = -\pdx\HH_1(\QQ^{ex}(t)) - \pdx\HH_2(\QQ^{im}(t)),
\end{align}
\end{subequations}
where $\QQ= (\qq,\vv_1,\vv_2)^T$, $\HH_1(\QQ^{ex}) = (\ff^\slow(\qq^{ex}), \bm 0, \bm 0)^T$ and $\HH_2 = (\vv_1^{im} + \vv_2^{im}, \AA_1^2 \qq^{im}, \AA_2^2 \qq^{im})^T$.
The superscripts \textit{ex,im} are added to highlight the explicit and implicit evaluation of the fluxes respectively.

As in the two-split case, system \eqref{eq:semi-implicit} is integrated using a scheme based on Butcher tableaux $((\tilde \alpha)_{r,k=1}^{s}, (\tilde \beta)_{r=1}^s, \cdot )$ for the explicit and $(\alpha_{r,k=1}^{s}, \beta_{r=1}^s, \cdot )$ for the implicit part.
At $t^n$ the variables $\QQ_{ex}$ and $\QQ_{im}$ are initialized by $\QQ^{ex,n} = \QQ^{im,n} = (\qq^n, \ff^\fast_1(\qq^n), \ff^\fast_2(\qq^n))^T$.
Then the $r$-th stage of the prediction for $r = 1, \dots, s$ in a time semi-discrete manner solving for the flux $\hh = \pdx \HH_1 + \pdx \HH_2$ is given by
\begin{subequations}\label{eq:semi-implicit-disc}
    \begin{eqnarray}
        &&\QQ^{ex,(r-1)} = \QQ^{ex,n} + \Delta t \sum_{k=1}^{r-1} \tilde{\alpha}_{rk} \hh^{(k)} \\
        &&\QQ^{im,(r-1)} = \QQ^{im,n} + \Delta t \sum_{k=1}^{r-1} \alpha_{rk} \hh^{(k)} \\
        &&\hh^{(r)} = \pdx\HH_1(\QQ^{ex,(r-1)}) + \pdx\HH_2(\QQ^{im,(r-1)} + \Delta t \alpha_{rr} \hh^{(r)}).
    \end{eqnarray}
\end{subequations}
Note that the fluxes in \eqref{eq:semi-implicit} for $\QQ^{ex}$ and $\QQ^{im}$ are the same and given by $\hh$ thus it constitutes a particular case of partitioned Runge-Kutta methods as detailed in \cite{Boscarino2016}.
Moreover, in \eqref{eq:semi-implicit-disc} the evolution of all variables including the relaxation variables $\vv_1, \vv_2$ are given.
To reduce the number of variables, the same procedure as in \eqref{eq:impl_sys_imex_3} can be applied yielding linear implicit systems with a Laplacian structure.
Moreover, at each stage, similar to \eqref{eq:impicit_2split_linsys_ho}, the projection $\vv_j^{(k)} = \ff_j^\fast(\qq^{(k)}), k = 1, \dots, r-1$ can be applied.
The scheme is completed by the correction update
\begin{equation}\label{eq:Relax_3Split_disc_update_ho}
    \ppsi^{n+1} = \ppsi^n - \frac{\Delta t}{\Delta x}\sum_{r=1}^{s} \left(\tilde \beta_{r}\left(\FF^{ex,(r)}_{i+1/2} - \FF^{ex,(r)}_{i-1/2}\right)+ \beta_r\sum_{j=1}^2\left(\FF^{im,(r)}_{j,i+1/2} - \FF^{im,(r)}_{j,i-1/2}\right)\right).
\end{equation}
We will use this procedure, together with a second order reconstruction of interface variables in space, to construct a second order three-split scheme for the MHD equations.
The details are given in Section \ref{sec:MHD}.
An assessment of the order of the numerical scheme has to be carefully addressed for each problem and flux splitting respectively and is beyond the scope of the paper.

\subsection{Multiple space dimensions}

For the derivation of the numerical scheme in two space dimensions, we resort to the first order two-split scheme for clarity.
All steps can be straightforwardly applied to the multi-split and higher order extensions from the previous sections.
We start from the two dimensional relaxation system \eqref{sys:RelaxMultiSplit_2D} with one slow and one fast sub-system associated to the directional fluxes $\ff^{1,\slow}, \ff^{2,\slow}, \ff^{1,\fast}, \ff^{2,\fast}$ in $x_1$- and $x_2$-direction respectively.
We consider a two dimensional computational domain $\Omega \subset \mathbb{R}^2$ which is paved by a Cartesian mesh of uniform grid size $\Delta x_l, l = 1,2 $ along each axis.
The cells are denoted by $C_I$ where $I = (i_1,i_2)$ on which the cell averages are defined by
\begin{equation}
    \qq_{i_1,i_2}(t) = \frac{1}{|C_I|}\int_{C_I} \qq(\xx,t) d\xx,
\end{equation}
where $\xx = (x_1,x_2)^T$.
The respective directional numerical fluxes will be denoted by $\FF^{1,im}, \FF^{2,im}$ defined by \eqref{eq:flux_im_center2} for the implicitly treated flux terms and $\FF^{1,ex}, \FF^{2,ex}$ defined by \eqref{eq:flux_ex_Rusanov} for the explicit ones.
Then the fully discrete scheme is constructed along the same lines as the 1D scheme by treating relaxation variables implicitly and slow dynamics explicitly.
This yields the following discretization
\begin{subequations}\label{sys:RelaxSplit_disc_2D}
    \begin{align}
      \begin{split}
          \ppsi^{(1)}_{i_1,i_2} - \ppsi^n_{i_1,i_2} &+ \frac{\Delta t}{\Delta x_1}(\FF^{1,ex,n}_{i_1+1/2,i_2} - \FF^{1,ex,n}_{i_1-1/2,i_2}) + \frac{\Delta t}{\Delta x_2}(\FF^{2,ex,n}_{i_1,i_2+1/2} - \FF^{2,ex,n}_{i_1,i_2-1/2}) \\
                & + \frac{\Delta t}{\Delta x_1}(\vv_{i_1+1/2,i_2}^{1,(1)} - \vv_{i_1-1/2,i_2}^{1,(1)}) + \frac{\Delta t}{\Delta x_2}(\vv_{i_1,i_2+1/2}^{2,(1)} - \vv_{i_1,i_2-1/2}^{2,(1)}) = 0, \label{eq:Relax_Split_disc_2D_psi}
      \end{split}\\
      \vv_{i_1+1/2,i_2}^{1,(1)} - \vv_{i_1+1/2,i_2}^{1,n} &+ \frac{\Delta t}{\Delta x_1}\AA_1^2 ( \ppsi_{i_1+1,i_2}^{(1)} - \ppsi_{i_1,i_2}^{(1)}) = 0,\label{eq:Relax_Split_disc_2D_v} \\
      \vv_{i_1,i_2+1/2}^{2,(1)} - \vv_{i_1,i_2+1/2}^{2,n} &+ \frac{\Delta t}{\Delta x_1}\AA_2^2 ( \ppsi_{i_1,i_2+1}^{(1)} - \ppsi_{i_1,i_2}^{(1)}) = 0,\label{eq:Relax_Split_disc_2D_w} \\
      \begin{split}
      \ppsi^{n+1}_{i_1,i_2} = \ppsi^n_{i_1,i_2} &- \frac{\Delta t}{\Delta x_1}(\FF^{1,ex,n}_{i_1+1/2,i_2} - \FF^{1,ex,n}_{i_1-1/2,i_2}) - \frac{\Delta t}{\Delta x_2}(\FF^{2,ex,n}_{i_1,i_2+1/2} - \FF^{2,ex,n}_{i_1,i_2-1/2}) \\
       &- \frac{\Delta t}{\Delta x_1} (\FF^{1,im,(1)}_{i_1+1/2,i_2} - \FF^{1,im,(1)}_{i_1-1/2,i_2}) - \frac{\Delta t}{\Delta x_2}(\FF^{2,im,(1)}_{i_1,i_2+1/2} - \FF^{2,im,(1)}_{i_1,i_2-1/2}). \label{eq:Relax_Split_disc_2D_update}
       \end{split}
    \end{align}
\end{subequations}
To reduce the number of variables, we replace $\vv^{1,(1)}$ and $\vv^{2,(1)}$ in \eqref{eq:Relax_Split_disc_2D_psi} by their respective updates given by \eqref{eq:Relax_Split_disc_2D_v} and \eqref{eq:Relax_Split_disc_2D_w}.
Further we make use of the relaxation equilibrium at $t^n$ which yields $\vv^{l,n} = \ff^{l,im}(\ppsi^n)$ for $l=1,2$.
Denoting the discrete weighted Laplacian by
\begin{equation}\label{eq:Laplacian}
    \LL(\ppsi_{i_1,i_2},\AA_1,\AA_2) = \frac{\AA_1^2}{\Delta x_1^2} \left(\ppsi_{i_1+1,i_2} - 2 \ppsi_{i_1,i_2} + \ppsi_{i_1-1,i_2}\right) + \frac{\AA_2^2}{\Delta x_1^2} \left(\ppsi_{i_1,i_2+1} - 2 \ppsi_{i_1,i_2} + \ppsi_{i_1,i_2-1}\right)
\end{equation}
the intermediate stage for $\ppsi^{(1)}$ writes
\begin{align*}
    \begin{split}
    \ppsi_{i_1,i_2}^{(1)} - \Delta t^2 \LL(\ppsi^{(1)}_{i_1,i_2}, \AA^2_1, \AA^2_2)= \ppsi^n_{i_1,i_2} &- \frac{\Delta t}{\Delta x_1}(\FF^{1,ex,n}_{i_1+1/2,i_2} - \FF^{1,ex,n}_{i_1-1/2,i_2}) - \frac{\Delta t}{\Delta x_2}(\FF^{2,ex,n}_{i_1,i_2+1/2} - \FF^{2,ex,n}_{i_1,i_2-1/2}) \\
    &- \frac{\Delta t}{\Delta x_1} (\FF^{1,im,n}_{i_1+1/2,i_2} - \FF^{1,im,n}_{i_1-1/2,i_2}) - \frac{\Delta t}{\Delta x_2}(\FF^{2,im,n}_{i_1,i_2+1/2} - \FF^{2,im,n}_{i_1,i_2-1/2}). \\
    \end{split}
\end{align*}
Along these lines also multi-dimensional high order schemes as given in the previous sections can be constructed.
In the following we will discuss two examples to illustrate the properties of the numerical scheme for an implicit-explicit flux splitting for the Euler equations and a flux splitting with two implicit fluxes for the MHD equations.

\section{Applications}

\subsection{The Euler equations}
To test the two-split numerical schemes, we consider the compressible Euler equations of gas dynamics.
The state vector is given by $\ppsi = (\rho, \rho \uu, \rho E)^T$, where $\uu = (u_1,\dots,u_d)^T \in \mathbb{R}^d$ denotes the velocity in a $d$ dimensional field, $\rho > 0$ the density and $E = e + E_\kin$ the specific total energy with $e> 0$ being the specific internal energy and $E_\kin = \frac{\uu \cdot \uu}{2}$ denotes the specific kinetic energy.
We apply the splitting of Toro \& V\'aszquez-Cend\'on \cite{Toro2012} which is given by
\begin{equation}\label{eq:Eulersplit}
    \ff^\slow = \begin{pmatrix}
        \rho \uu \\ \rho \uu \otimes \uu \\ \rho E_\kin \uu
    \end{pmatrix}, \quad
    \ff^\fast = \begin{pmatrix}
        0 \\ p \\ (\rho e + p ) \uu
    \end{pmatrix}.
\end{equation}
To close the system we have to define an equation of state which connects the pressure to the density and internal energy.
The framework of the numerical scheme would also allow the usage of a general equation of state, see e.g. \cite{Thomann2023,Boscheri2021}.
Here, we consider for simplicity an ideal gas law, where the pressure is given by
\begin{equation}
    p (\rho, e)= (\gamma - 1) \rho e, \quad \gamma > 1.
\end{equation}

The eigenvalues in normal flow direction for the convective sub-system associated to $\ff^\slow$ are given by
\begin{equation}\label{eq:EV_ex}
    \lambda^\slow = \uu \cdot \nn
\end{equation}
which has multiplicity $2 + d$.
The eigenvalues in normal direction with respect to the pressure sub-system associated to $\ff^\fast$ are given by
\begin{equation}\label{eq:EV_p}
    \lambda_1^\fast = 0, \quad \lambda_{2,3}^\fast = \frac{1}{2}(\uu \cdot \nn \pm \sqrt{(\uu \cdot \nn)^2 + 4 c^2}),
\end{equation}
where $c$ denotes the sound speed given by $c^2 = \gamma\frac{p}{\rho}$ for an ideal gas.
Note that the eigenvalues of the convective sub-system are independent of the pressure, i.e. from the Mach number which is given as the ratio between local flow speed $|\uu|$ and the sound speed $c$ which itself depends on the pressure.
Thus the CFL condition \eqref{eq:CFL} depends solely on $\uu$ independently of the Mach number regime.

Note that different splitting approaches for the full Euler equations were proposed, e.g. \cite{BoscarinoRussoScandurra2018,Bispen2017} and references therein.
We have chosen the splitting from \cite{Toro2012} as it separates the pressure dependent flux terms completely from the velocity based convective terms.
In particular, for constant pressure and velocity, the pressure sub-system is constant, thus variations in the density are evolved entirely by the explicit material sub-system, see \cite{Boscheri2020,Boscheri2021}.

Another property which is important in the context of flows in different Mach number regimes is the so-called asymptotic preserving property.
This means, as the Mach number tends to zero, the numerical scheme should yield a consistent discretization of the limit system given by the incompressible Euler equations, see e.g. \cite{Dellacherie2010,KlainermanMajda1981,Klein1995}.

We start with the contact preserving property.

\begin{lemma}[Contact preserving property]\label{lem:contact}
    Let $u^0 = u_0$ constant, as well as $p^0 = p_0$.
    Then, for the first order scheme \eqref{eq:impl_sys_imex}, \eqref{eq:Relax_Split_disc_update}, the velocity and pressure is preserved, i.e. $u^n = u_0$, $p^n = p_0$ for any time $t^n = n \Delta t$ under the CFL condition \eqref{eq:CFL}.
\end{lemma}
\begin{proof}
    For simplicitly we consider the one-dimensional framework where $d=1$, i.e. $\uu = u \in \mathbb{R}$. However the analogous steps can be performed in the two or higher dimensional setting.

    Let us consider one time step from $t = t^0$ to $t^1$ given by $\Delta t$ under the CFL condition \eqref{eq:CFL} where the initial data is given by $\qq^0 = (\rho^0, \rho^0 u_0, e^0(\rho^0,p_0) + \frac{1}{2}\rho^0 u_0^2)^T$.
    Let $\AA_\fast = a_\fast \mathbb{I}$.
     For $u_0 = 0$ we immediately find $\qq^{n+1} = \qq^n$. Let further w.l.o.g. $u_0 > 0$ constant.
    Then the update of the density reads for the first order scheme
    \begin{equation}\label{eq:Contact_density_update}
        \rho^1_i = \rho^0_i + \frac{\Delta t}{\Delta x}u_0 (\rho_{i}^0 -\rho_{i-1}^0).
    \end{equation}
    For the intermediate states, we have for the density
    \begin{equation}\label{eq:Lemma_aux_1}
        \rho^{(1)}_i - \frac{\Delta t^2}{\Delta x^2} a_\fast^2 \left(\rho ^{(1)}_{i+1} - 2 \rho_i^{(1)} + \rho_{i-1}^{(1)}\right)= \rho^0_i - u_0\frac{\Delta t}{\Delta x}(\rho_i^0 - \rho_{i-1}^0)
    \end{equation}
    for the momentum
    \begin{align}
        (\rho u)^{(1)}_i - \frac{\Delta t^2}{\Delta x^2} a_\fast^2 \left((\rho u)^{(1)}_{i+1} - 2 (\rho u)_i^{(1)} + (\rho u)_{i-1}^{(1)}\right)&= u_0 \left(\rho^0_i - u_0\frac{\Delta t}{\Delta x}(\rho_i^0 - \rho_{i-1}^0)\right)\\
        &= \rho^{(1)}_i u_0 - \frac{\Delta t^2}{\Delta x^2} a_\fast^2 u_0 \left(\rho^{(1)}_{i+1} - 2 \rho_i^{(1)} + \rho_{i-1}^{(1)}\right),
    \end{align}
    i.e.
    \begin{align}
        (\rho (u - u_0))^{(1)}_i - \frac{\Delta t^2}{\Delta x^2} a_\fast^2 \left((\rho (u - u_0))^{(1)}_{i+1} - 2 (\rho (u - u_0))_i^{(1)} + (\rho (u - u_0))_{i-1}^{(1)}\right)
        = 0 \quad \forall i = 1, \dots, N.
    \end{align}
     This relation can be rewritten as $\CC [\rho^{(1)}_i (u^{(1)}_i - u_0)]_{i=1}^N = 0$ with a symmetric tridiagonal positive definite matrix $\CC\in \mathbb{R}^{N\times N}$, i.e. $\CC$ is invertible.
     Thus either $\rho_i^{(1)} = 0$ or $u^{(1)}_i = u_0$ for all $i$.

     To rule out the first case, we show that the density indeed stays positive for positive $\rho^0_i$.
     The right hand side of \eqref{eq:Lemma_aux_1} can be written as convex combination of positive states $\rho_i^0, \rho_{i-1}^0$ since due to the CFL condition $\Delta t/\Delta x \leq \nu/u_0$
     \begin{equation}
        \rho^0_i - u_0\frac{\Delta t}{\Delta x}(\rho_i^0 - \rho_{i-1}^0) \geq (1-\nu) \rho_i^n + \rho_{i-1}^n > 0,
     \end{equation}
     for all $i$ and $\nu$ from the CFL condition \eqref{eq:CFL}.
     Thus, $\CC[\rho_i^{(1)}]_{i=1}^N > 0$. Since $\CC^{-1}$ is symmetric positive definite for the intermediate density holds $\rho^{(1)}_i > 0$ for all $i$.
     Thus $u_i^{(1)}= u_0$.
     For the total energy, we have with $E_\kin^0 = \frac{u_0^2}{2}$ the following implicit relation
     \begin{align}
         (\rho E)^{(1)}_i - \frac{\Delta t^2}{\Delta x^2} a_\fast^2 \left((\rho E)^{(1)}_{i+1} - 2 (\rho E)_i^{(1)} + (\rho E)_{i-1}^{(1)}\right)=  E_\kin^0 \left(\rho^0_i - u_0\frac{\Delta t}{\Delta x}(\rho_i^0 - \rho_{i-1}^0)\right) + \frac{p_0}{\gamma - 1} \\
         = (\rho)^{(1)}_i E_\kin^0 + \frac{p_0}{\gamma - 1} - \frac{\Delta t^2}{\Delta x^2} a_{\fast}^2 \left((\rho )^{(1)}_{i+1} E_\kin^0 - 2 (\rho)_i^{(1)} E_\kin^0 + (\rho)_{i-1}^{(1)} E_\kin^0\right).
     \end{align}
     Since $(\rho E)^{(1)} = \rho^{(1)}E_\kin^0 + \frac{p^{(1)}}{\gamma - 1}$ and $p_0$ constant, it simplifies to a relation for the pressure as follows
     \begin{align}
         (p-p_0)^{(1)}_i - \frac{\Delta t^2}{\Delta x^2} a_{\fast}^2 \left((p-p_0)^{(1)}_{i+1} - 2 (p-p_0)_i^{(1)} + (p-p_0)_{i-1}^{(1)}\right)= 0
     \end{align}
     i.e. $\CC[(p-p_0)^{(1)}_i]_{i=1}^N = 0$.
     With the analogous argumentation as detailed above, we find $p^{(1)} = p_0$.
     For the correction it holds according to the final numerical update \eqref{eq:Relax_Split_disc_update} that the implicit flux contributions vanish, and we have the density update \eqref{eq:Contact_density_update} and
     \begin{equation}
         \rho^1 u^1 = u_0 \rho^1
     \end{equation}
     i.e. $u^1 = u_0$ and
     \begin{equation}
         p^1 + (\gamma - 1) \rho^1 E_\kin^1 = p_0 + (\gamma - 1)\rho^1 E_\kin^0,
     \end{equation}
    i.e. $p^1 = p^0$.
    Thus analogously it holds passing from $t^{n-1}$ to $t^{n}$ for $n> 1$ which concludes the proof.
\end{proof}
Note, that we have assumed in the proof, that the linear systems are solved exactly. In practise some numerical errors can accumulate, i.e. the velocity and pressure errors are small but might be larger than machine precision.

We turn to the asymptotic consistency with the incompressible Euler equations.
Therefore, we define the Mach number $M = |\uu|/c$ which by rescaling arguments appears in the Euler equations as follows
\begin{equation}\label{sys:EulerM}
    \displaystyle
    \pdt
    \begin{pmatrix}
        \rho\\\rho\uu\\\rho E
    \end{pmatrix}
    +
    \divv
    \begin{pmatrix}
        \rho \uu \\ \rho \uu \otimes \uu \\ M^2\rho E_\kin \uu
    \end{pmatrix}
    + \divv \begin{pmatrix}\displaystyle
        0 \\ \displaystyle\frac{p}{M^2} \\ {(\rho e + p )} \uu
    \end{pmatrix} = 0
\end{equation}
where a rescaled total energy is given by $\rho E = M^2\rho E_\kin + {\rho e}$.
For this system, we consider so-called well-prepared initial data, i.e. let the initial condition of the continuous equations be given by
\begin{align}\label{def:wpdata}
    \rho = \rho_0 + \mathcal{O}(M^2), \quad \rho_0 \text{ const,} \quad \uu = \uu_0 + \mathcal{O}(M^2),\quad \nabla \cdot \uu_0 = 0, \quad p = p_0 + \mathcal{O}(M^2), \quad p_0 \text{ const.}
\end{align}
for which was formally shown in e.g. \cite{Dellacherie2010} that as $M\to 0$, the system formally converges to the incompressible Euler equations
\begin{equation}
   \pdt \uu_0 + \uu_0 \cdot \nabla \uu_0 + \frac{\nabla p_2}{\rho_0} = 0, \quad \nabla \cdot \uu_0 = 0.
\end{equation}
Note that this formulation for sufficiently smooth solutions is equivalent to solving
\begin{equation}\label{def:wpdata2}
    \divv(\uu_0 \cdot \nabla \uu_0) = \frac{1}{\rho_0}\Delta p_2.
\end{equation}
Thus discretized well-prepared data should fulfil \eqref{def:wpdata} and \eqref{def:wpdata2} also on the discrete level.
Since this is a truly multi-dimensional problem, we consider the two-dimensional set-up.
\begin{lemma}
    Consider discretized well-prepared initial data in the sense of \eqref{def:wpdata} and \eqref{def:wpdata2} up to $\mathcal{O}(h)$ where $h$ denotes a uniform space discretization on Cartesian grids.
    Then the first order scheme given by \eqref{eq:impl_sys_imex}, \eqref{eq:Relax_Split_disc_update} is asymptotic preserving in the sense that for $M\to 0$ the discretization yields a consistent discretization of the incompressible Euler equations at first order.

    Moreover, the solution obtained after one iteration of the scheme is again discretized well-prepared data up to $\mathcal{O}(h)$.
\end{lemma}
\begin{proof}
    Let the initial data be given by \eqref{def:wpdata} and let us consider one time step from $t=t^n$ to $t^{n+1}$.
    We consider the numerical scheme in two dimensions on a Cartesian grid with homogeneous step size $\Delta x_l = h$, $l=1,2$ and $\AA_1 = \AA_2 = a_\fast \mathbb{I}$.
    Further, let here $\uu = (u^1,u^2)^T$.
    Let $L$ denote the discrete Laplacian operator
    \begin{equation}
        L(\cdot)_{i_1,i_2} = \frac{1}{h^2}\left((\cdot)_{i_1+1,i_2} - 2 (\cdot)_{i_1,i_2} + (\cdot)_{i_1-1,i_2} + (\cdot)_{i_1,i_2+1} - 2 (\cdot)_{i_1,i_2} + (\cdot)_{i_1,i_2-1}\right),
    \end{equation}
    a discrete divergence operator as
    \begin{equation}
        \divv(\uu_{i_1,i_2}) = \frac{1}{2h} \left(u_{i_1+1,i_2}^1 - u_{i_1-1,i_2}^1 + u_{i_1,i_2+1}^2 - u_{i_1,i_2-1}^2\right)
    \end{equation}
    and diffusion operator as
    \begin{equation}
        D(\cdot)_{i_1,i_2} = \frac{h}{2}\max_{i_1,i_2}|\uu_{i_1,i_2}|L(\cdot)_{i_1,i_2})
    \end{equation}
    Then, the first stage of the numerical scheme in two dimensions for the scaled system \eqref{sys:EulerM} with global diffusion in the explicit fluxes is given by
    \begin{align*}
            \rho^{(1)}_{i_1,i_2} - \Delta t^2 a_\fast^2 L\rho^{(1)}_{i_1,i_2} &= \rho^n_i - \Delta t \divv(\rho \uu)_{i_1,i_2}^n + \Delta t D\rho^n_{i_1,i_2} \\
            (\rho u^1)^{(1)}_{i_1,i_2} - \Delta t^2 a_\fast^2 L(\rho u^1)^{(1)}_{i_1,i_2} &= (\rho u^1)^n_{i_1,i_2} - \Delta t\divv(\rho \uu u^1)^n_{i_1,i_2} - \frac{\Delta t}{2h M^2} (p_{i_1+1,i_2}^n - p_{i_1-1,i_2}^n)+ \Delta t D(\rho u^1)_{i_1,i_2}^n\\
            (\rho u^2)^{(1)}_{i_1,i_2} - \Delta t^2 a_\fast^2 L(\rho u^2)^{(1)}_{i_1,i_2} &= (\rho u^2)^n_{i_1,i_2} - \Delta t \divv(\rho \uu u^2)^n_{i_1,i_2} - \frac{\Delta t}{2h M^2} (p_{i_1,i_2+1}^n - p_{i_1,i_2-1}^n) + \Delta tD(\rho u^2)^n_{i_1,i_2}\\
            (\rho E)^{(1)}_{i_1,i_2} - \Delta t^2 a_\fast^2 L(\rho E)^{(1)}_{i_1,i_2} &= (\rho E)^n_{i_1,i_2} - M^2\Delta t \divv(\rho E_\kin \uu)^n_{i_1,i_2} - {\Delta t}\divv(\frac{\gamma}{\gamma - 1} p^n \uu^n)_{i_1,i_2} + \Delta t D(\rho E)^n_{i_1,i_2}.
    \end{align*}
    Making use of the definition of well-prepared initial data and assume that $\uu_0$ at $t^n$ fulfils the discrete divergence operator up to $\mathcal{O}(h)$, we obtain
        \begin{align}
        &\rho^{(1)}_{i_1,i_2} - \Delta t^2 a_\fast^2 L\rho^{(1)}_{i_1,i_2} = \rho_{0} + \mathcal{O}(h\Delta t,M^2), \\
        &\left(\frac{p}{\gamma - 1} + M^2\rho E_\kin\right)^{(1)}_{i_1,i_2} - \Delta t^2 a_\fast^2 L\left(\frac{p}{\gamma - 1}+ M^2\rho E_\kin\right)^{(1)}_{i_1,i_2} = \frac{p_{0}}{\gamma - 1} + \mathcal{O}(h\Delta t, M^2).
    \end{align}
    Since $a_\fast$ is an approximation of the largest absolute eigenvalue of the pressure sub-system, it scales with $1/M$.
    Thus $a^2_{\fast} = \mathcal{O}(M^{-2})$ and the inverse of the coefficient matrix $(I - \Delta t^2 a_\fast^2 L)$ has eigenvalues that scale either of order $\mathcal{O}(1)$ or $\mathcal{O}(M^2)$ with respect to the Mach number for $h/\Delta t = \mathcal{O}(1)$, see a similar technique done in \cite{Bispen2017}.
    Thus the density expansion is given by $\rho^{(1)} = \rho_0 + \mathcal{O}(h\Delta t, M^2)$, $p^{(1)} = p_0 + \mathcal{O}(h\Delta t, M^2)$.
    Regarding the velocity terms, we have
    \begin{align*}
        (\rho u^1)^{(1)}_{i_1,i_2} - \Delta t^2 a_\fast^2 L(\rho u^1)^{(1)}_{i_1,i_2} = \rho_0( u_0^1)^n_{i_1,i_2}& - \Delta t\rho_0 \divv(\uu_0 u_0^1)^n_{i_1,i_2} \\&- \frac{\Delta t}{2h} (p_{2,i_1+1,i_2}^n - p_{2,i_1-1,i_2}^n) + \rho_0 \Delta t D(u_0^1)^n_{i_1,i_2}+ \mathcal{O}(M^2), \\
        (\rho u^2)^{(1)}_{i_1,i_2} - \Delta t^2 a_\fast^2 L(\rho u^2)^{(1)}_{i_1,i_2} = \rho_0( u_0^2)^n_{i_1,i_2} &- \Delta t \rho_0\divv( \uu_0 u_0^2)^n_{i_1,i_2} \\
        &- \frac{\Delta t}{2h} (p_{2,i_1,i_2+1}^n - p_{2,i_1,i_2-1}^n) + \rho_0 \Delta t D (u_0^2)^n_{i_1,i_2}+ \mathcal{O}(M^2).
    \end{align*}
    With $\divv(\uu_0)^n_{i_1,i_2} = \mathcal{O}(h)$ it follows for a scalar quantity $\phi$
    \begin{align*}
        \divv(\uu_0 \phi)_{i_1,i_2} &= \left(\frac{\phi_{i_1+1,i_2} - \phi_{i-1,i_2}}{2h}\right)\frac{u_{0,i_1+1,i_2}^1 + u_{0,i_1-1,i_2}^1}{2} + \left(\frac{\phi_{i_1,i_2+1} - \phi_{i_1,i_2-1}}{2h}\right)\frac{u_{0,i,i_2+1}^2 + u_{0,i,i_2-1}^2}{2} \\
        &+ \frac{\phi_{i_1+1,i_2} - 2 \phi_{i_1,i_2} + \phi_{i_1-1,i_2}}{h^2}\left(\frac{h^2}{2}\frac{u_{0,i_1+1,i_2}^1 - u_{0,i_1-1,i_2}^1}{2h}\right)\\
        &+ \frac{\phi_{i_1,i_2+1} - 2 \phi_{i_1,i_2} + \phi_{i_1,i_2-1}}{h^2}\left(\frac{h^2}{2}\frac{u_{0,i_1,i_2+1}^2 - u_{0,i_1,i_2-1}^2}{2h}\right) + \mathcal{O}(h)
    \end{align*}
    which is consistent with $\uu_0 \cdot \nabla \phi + h^2\pdxl{1}( u \pdxl{1} \phi) + h^2\pdxl{2}( v \pdxl{2} \phi) + \mathcal{O}(h) = \uu_0 \cdot \nabla \phi + h^2 \nabla \cdot (\uu_0 \cdot \nabla \phi) + \mathcal{O}(h) = \uu_0 \cdot \nabla \phi + \mathcal{O}(h)$.
    Note that the diffusion term $D(\phi)_{i_1,i_2}$ is of order $\mathcal{O}(h)$ for bounded velocity.
    Thus with abuse of notation, we can write
    \begin{align}
    (\rho u^1)^{(1)}_{i_1,i_2} - \Delta t^2 a_\fast^2 L(\rho u^1)^{(1)}_{i_1,i_2}  = \rho_0( u_0^1)^n_{i_1,i_2} &- \Delta t\rho_0 \overline{\uu_0}^n_{i_1,i_2} \cdot \nabla_{i_1,i_2} u_0^{1,n} \\
    &- \frac{\Delta t}{2h} (p_{2,i_1+1,i_2}^n - p_{2,i_1-1,i_2}^n) + \mathcal{O}(h \Delta t) + \mathcal{O}(M), \\
    (\rho u^2)^{(1)}_{i_1,i_2} - \Delta t^2 a_\fast^2 L(\rho u^2)^{(1)}_{i_1,i_2} = \rho_0( u_0^2)^n_{i_1,i_2} &- \Delta t\rho_0 \overline{\uu_0}^n_{i_1,i_2} \cdot \nabla_{i_1,i_2} u_0^{2,n}\\ &- \frac{\Delta t}{2h} (p_{2,i_1,i_2+1}^n - p_{2,i_1,i_2-1}^n) + \mathcal{O}(h\Delta t)+ \mathcal{O}(M).
    \end{align}
    With the same argumentation as above, we have
    as $M \to 0$
    \begin{equation}\label{eq:intermediate_veloc_AP}
        \uu^{(1)}_{0,i_1,i_2} = \uu_{0,i_1,i_2}^n - \Delta t \overline{\uu_0}_{i_1,i_2}^n \cdot \nabla_{i_1,i_2}\uu_0^n + \frac{\Delta t}{\rho_0}\nabla_{i_1,i_2} p_{2}^n \mathbb{I} + \mathcal{O}(h\Delta t)
    \end{equation}
    which is a consistent first order discretization of the incompressible Euler equations.

    Assuming a CFL condition which $\Delta t/h = \mathcal{O}(1)$ we have $h\Delta t = \mathcal{O}(h^2)$.
    Regarding the divergence free property of the velocity field $\uu^{(1)}_{i_1,i_2}$, we apply the discrete divergence on \eqref{eq:intermediate_veloc_AP}.
    Since the data at $t^n$ is well-prepared up to $\mathcal{O}(h)$, we have
    \begin{equation}\label{eq:disc_wb_constr_u}
        \divv(- \overline{\uu_0}^n \cdot \nabla_{i_1,i_2}\uu_0^n + \frac{1}{\rho_0}\nabla_{i_1,i_2} p_{2}^n)_{i_1,i_2} = \mathcal{O}(h).
    \end{equation}
    Thus we obtain
    \begin{equation}
        \divv(\uu_{0}^{(1)})_{i_1,i_2} = \mathcal{O}(h).
    \end{equation}
    Thus the solution after the prediction step is well-prepared up to first order discretization errors $\mathcal{O}(h)$.

    Then, in the correction step, we have using the predicted variables $\qq^{(1)}_{i_1,i_2}$ derived above
    \begin{align}\label{eq:disc_update_limit}
        \rho^{n+1}_{i_1,i_2} &= \rho_0 + \mathcal{O}(h\Delta t, M^2), \\
        (\uu_{0})_{i_1,i_2}^{n+1} &= (\uu_{0})^n_{i_1,i_2} - \Delta t \overline{\uu_0}_{i_1,i_2}^n \cdot \nabla_{i_1,i_2}\uu_0^n + \frac{\Delta t}{\rho_0}\nabla_{i_1,i_2} p_{2}^{(1)} \mathbb{I} + \mathcal{O}(h\Delta t,M), \label{eq:disc_update_limit_u}\\
        p^{n+1} & = p_0 + \mathcal{O}(h\Delta t, M^2).
    \end{align}
    Thus as the Mach number tends to zero, a consistent first order discretization of the incompressible Euler equations is obtained.
    Finally we have, applying the discrete divergence operator on \eqref{eq:disc_update_limit_u} and using \eqref{eq:disc_wb_constr_u},  $\divv(\uu_0^{n+1})_{i_1,i_2} = \mathcal{O}(h)$.
    Thus overall we have a first order approximation of the well-prepared initial data at $t^{n+1}$ and a consistent discretization with the limit model.
    This concludes the proof.
\end{proof}

\subsection{The ideal MHD equations}\label{sec:MHD}
To test the three-split numerical schemes, we consider the ideal MHD equations which exhibit formally three scales, the material scale associated with $\uu$, a pressure sub-system associated with the Mach number $M$ and a magnetic sub-system associated to the strength of the magnetic field with respect to the local fluid velocity characterised by the Alfv\'en number denoted by $M_A$.

Compared to the Euler equations considered above, an additional difficulty is posed by the solenoidal property of the magnetic field, i.e. the divergence error of the magnetic field has to be controlled if not preserved exactly by the numerical scheme.
Since in this work, we are mainly interested in providing a general framework, we do not apply specific divergence free operators as done in the recent works \cite{Boscheri2024} or \cite{Fambri2021}.
In contrast, we make use of the divergence cleaning approach by Munz \cite{Munz2000}, which means that we add an additional variable $\phi$ whose task it is to propagate the divergence errors out of the computational domain with a constant cleaning speed $c_h^2$ and at the same time dampen the divergence errors.
Thus, the approach fits naturally into the relaxation framework utilized here to construct the numerical schemes.

We apply the MHD splitting proposed by Fambri in \cite{Fambri2021} which can be viewed as an extension of the Toro-V\'azquez-Cend\'on splitting of the Euler equations.
An additional sub-system appears which contains all the dynamics of the magnetic field denoted by $\BB = (B_1, \dots, B_d)^T \in \mathbb{R}^{d}$.
Including the variable $\phi$ for the divergence cleaning, the state vector is given by $\ppsi = (\rho, \rho\uu, \rho E, \BB, \phi)^T $ and the three fluxes read
\begin{equation}\label{eq:MHDsplit}
    \ff^\slow = \begin{pmatrix}
        \rho \uu \\ \rho \uu \otimes \uu \\ \rho E_\kin \uu \\\bm 0 \\0
    \end{pmatrix}, \quad
    \ff^\fast_1 = \begin{pmatrix}
        0 \\ p \\ (\rho e + p ) \uu \\ \bm{0} \\ 0
    \end{pmatrix}, \quad
    \ff^\fast_2 = \begin{pmatrix}
    0 \\ m \mathbb{I} - \frac{1}{4\pi} \BB \otimes \BB \\ 2 m \uu - \frac{1}{4\pi} \BB (\uu \cdot \BB)\\ \BB \otimes \uu - \uu \otimes \BB + c_h^2 \phi \mathbb{I} \\ \BB
    \end{pmatrix}
\end{equation}
where now $\rho E = \rho E_\kin + \rho e + m$ with the magnetic pressure contribution $m = \frac{\|\BB\|^2}{8\pi}$.
Therein $e$ and $E_\kin$ are defined as in the Euler equations using and ideal gas law.
The eigenvalues for the material sub-system associated to the flux $\ff^\slow$ and the ones for the pressure sub-system associated to $\ff_1^\fast$ are identical up to additional $0$ eigenvalues with the ones of the respective sub-systems of the Euler equations and are given by \eqref{eq:EV_ex} and \eqref{eq:EV_p} respectively.
The characteristic speeds associated with the magnetic sub-system, i.e. eigenvalues of $\nabla_\qq \ff^{\fast}_2$ are given by
\begin{align}
    &\lambda_{2,8}^{2,\fast} = \frac{1}{2} \left(\uu \cdot \nn \mp \sqrt{(\uu \cdot \nn)^2 + 4 b^2}\right), \quad \lambda_{3,7}^{2,\fast} = \frac{1}{2} \left(\uu \cdot \nn \mp \sqrt{(\uu \cdot \nn)^2 + 4 b_n^2}\right),\\
    &\lambda_{1,9}^{2,\fast} = \pm c_h, \quad \lambda_{4,5,6}^{2,\fast} = 0,
\end{align}
where $b = \|\BB\|^2/\sqrt{4\pi\rho}$ and $b_n = (\BB \cdot \nn)/\sqrt{4\pi\rho}$ denote the total Alfv\'en speed and the directional Alfv\'en speed respectively.
The Alfv\'en number is then defined as $M_A = |\uu|/b$.
Note that the characteristic speeds of the magnetic sub-system do not depend on the Mach number and the characteristic speeds of the pressure sub-system do not depend on the Alfv\'en number. Thus this splitting is separating the three scales of the MHD system.
This is also the case for the slightly different splitting considered in \cite{Boscheri2024}.
Therefore, the CFL condition, as in the Euler case, only depends on the material waves and is given by \eqref{eq:CFL}.

The choice of the constant cleaning speed $c_h$ is not unique.
Here, we choose it larger than the maximal absolute eigenvalue of the magnetic sub-system.
However, note, that, since the diffusion of the numerical scheme for the magnetic sub-system is of the order of the maximal absolute eigenvalue $\lambda^{2,\fast}$, choosing it arbitrarily large implies an increased diffusion on the magnetic field.
In practise, it is usually taken twice the $\max |\lambda_{2,8}^{2,\fast}|$.
We set $\AA_1 = a_1^\fast \mathbb{I}$ and $\AA_2 = a_2^\fast \mathbb{I}$ where $a_1^\fast$ and $a_2^\fast$ satisfy the sub-characteristic conditions from Lemma \ref{lem:subchar_multisplit_1D} in the one-dimensional case and \eqref{eq:subchar_multisplit_2D} in the two-dimensional case.

Since the magnetic field is entirely governed by the subsystem connected to $\ff^\fast_2$, the magnetic subsystem will be resolved before the pressure subsystem, i.e. first the update with $\ff^\fast_2$ in \eqref{eq:impl_sys_imex_3_s1star} and then $\ff^\fast_1$ is performed in \eqref{eq:impl_sys_imex_3_s1}, similar to the three-split scheme derived in \cite{Boscheri2024}.
Then the predicted update of the magnetic field $\BB^{(1)} = \BB^{(1,\ast)}$ can be used to update the magentic field contributions in the momentum and total energy equation in the stage \eqref{eq:impl_sys_imex_3_s1}.
This means that due to the particular splitting applied here, only $d$ instead of in general $2d + 3$ implicit systems have to be solved in \eqref{eq:impl_sys_imex_3_s1star}.

Since the Fambri-splitting of the MHD system \eqref{eq:MHDsplit} is an extension of the Toro-V\'azquez-Cend\'on-splitting \eqref{eq:Eulersplit}, we recover the contact preserving property for constant pressure, velocity and magnetic field.
The proof is analogous to the one presented for the Euler equations and is thus omitted.
For a fixed Alfv\'en number regime, the asymptotic consistency with the incompressible Euler equations in a magnetic field as $M\to 0$ is recovered.
The analogous result on the asymptotic preserving property with respect to the Mach number is obtained.
The proof is analogous to the one presented above and is thus omitted.

\section{Numerical results}
\label{sec:NumRes}
In the following we present a suite of test cases for the Euler equation and the ideal MHD equations to numerically asses accuracy and robustness of the two- and three-split numerical scheme at a wide range of Mach and Alfv\'en numbers.
If not otherwise specified, we set $\gamma = 1.4$ and the time step is computed with the CFL condition \eqref{eq:CFL} with a material CFL number $\nu$.
To allow a comparison to the CFL condition of a fully explicit scheme, we define the CFL condition
\begin{equation}
    \Delta t = \nu_{ex} \frac{\Delta x}{\max_i |\lambda(\qq_i^n)|}
\end{equation}
with explicit CFL number $\nu_{ex}$.
If not otherwise mentioned, all test cases are performed with the second order scheme employing the LSDIRK2 scheme from \cite{AvgBerIolRus2019} which in Butcher notation is given by
\begin{equation}
    \renewcommand{\arraystretch}{1.25}
    \begin{array}{c|cc}
        0 & 0 & 0\\
        \tilde{c} & \tilde{c} & 0 \\\hline
        & 1- \alpha & \alpha
    \end{array} \quad
 \begin{array}{c|cc}
    \eta & \eta & 0\\
    1 & 1-\eta & \eta \\\hline
    & 1- \eta & \eta
\end{array},
\end{equation}
with $\eta = 1- 1/\sqrt{2}$ and $\tilde c = 1/(2\eta)$.
To obtain overall second order accuracy, a TVD reconstruction is carried out for a quantity $q(x)$ in each spatial direction, hence obtaining a piecewise polynomial representation $r(x)$ which, in $x$-direction, writes
\begin{equation}
    r_i(x) = c_0 q_i + c_1 (x - x_i),
\end{equation}
with the expansion coefficients computed using a minmod limiter \cite{Toro2009} as
\begin{equation}
    c_0 = 1, \quad c_1 = \text{minmod}\left(\frac{q_{i+1}-q_i}{\Delta x}, \frac{q_i - q_{i-1}}{\Delta t}\right).
\end{equation}
A second order numerical flux is then obtained setting $q_{i+1/2}^+ = r_{i+1}(x_{i+1/2})$ and $q_{i+1/2}^- = r_i(x_{i-1/2})$ as arguments of the Rusanov flux \eqref{eq:flux_ex_Rusanov}.
The resulting scheme is labelled \IMEX which is short for Semi-Implicit Finite Volume Prediction-Correction scheme.

\subsection{The Euler equations}
To complete the two-split scheme for the Euler equations, we have to define the matrix $\AA^\fast$ which for all test cases is given by $\AA^\fast = a \mathbb{I}$ with $a = \max_{i}\max_j|\lambda^\fast_j(\qq_i^n)|$ where $i$ runs over all cell averages and $j$ all characteristic speeds.

\subsubsection{Travelling isentropic vortex}
An exact solution of the Euler equations is given by the travelling isentropic vortex in two space dimensions, i.e. $\uu = (u_1,u_2)^T$.
To verify the convergence of the second-order two-split scheme of the Euler equations, we consider the isentropic vortex scaled with respect to the Mach number from \cite{MichelDansac2022} shifted by the constant velocity $\uu_m$.
Let $\xx = (x_1,x_2)^T$. Then the initial condition
\begin{equation}
    \label{eq:vortex}
    \begin{cases}\displaystyle
        \rho^0(\xx) = 1 - \frac{M^2}{8} \exp({- 2 a^2 \|\xx\|^2}), \quad p^0(\xx) = \frac{1}{M^2}\rho(\xx)^\gamma,                                      \\[5mm]\displaystyle
        \uu^0(\xx) = \uu_m + a \sqrt{\frac \gamma 2} \exp({- a^2 \|\xx\|^2}) \rho^0(\xx)^{\frac{\gamma}{2} - 1} \begin{pmatrix} x_2 \\ -x_1 \end{pmatrix}.
    \end{cases}
\end{equation}
with leads to the exact solution
\begin{equation}
    \rho(\xx,t) = \rho^0(\xx - \uu_m t), \quad \uu(\xx,t) = \uu^0(\xx - \uu_m t), \quad p(\xx,t) = p^0(\xx - \uu_m t).
\end{equation}
For the simulation, we set $a = 8$ and $\uu_m = (1,1)^T$ and consider periodic boundary conditions.
We choose the computational domain $[-1,1]\times[-1,1]$.
For the assessment of the convergence, we choose $M=10^{-k}, k = 1,\dots,3$ and moderate CFL numbers, i.e. $\nu_{mat}/\nu_{ex} = 10$ and $\nu_{mat} = 2.5 M$.
We calculate the same number of time steps for all Mach numbers, i.e. the final time is given by $t_f = M$.
The experimental order of convergence of the second order two-split scheme for this test case is as expected, see the errors and convergence rates reported in Table \ref{tab:EOCEuler}.
Note that the vortex solution has the following expansion with respect to the Mach numer $\rho^0 = \rho_0 + M^2\rho_2, \nabla\cdot \uu^0 = 0$ and $p^0 = p_0 + M^2p_2 + \mathcal{O}(M^4)$ in the sense of the scaled equations \eqref{sys:EulerM}.
We rerun the test case with $\uu_m = (0,0)^T$ at fixed material CFL number $\nu_{mat} = 0.1$ for Mach numbers $M = 10^{-k}, k=0,\dots,4$ and calculate the relative error of the density and pressure with respect to the leading order terms $\rho_0$ and $p_0$ at time $t_f = 0.1$.
The results are reported in Table \ref{tab:EulerAP} which verifies that the error in density and pressure decrease with $\mathcal{O}(M^2)$ as the Mach number decreases.

\begin{table}[htpb]
    \centering
    \renewcommand{\arraystretch}{1.25}
    \begin{tabular}{cccccccc}
        $M$ & $N$ &  $|\uu|$ & & $\rho $ & & $p$ & \\ \hline
        \multirow{4}{*}{$10^{-1}$} & 32	& 8.75E-03 & ---   &1.56E-03 &   ---  &	9.73E-03 &	      ---             \\
        & 64	& 2.71E-03 & 1.688	&6.26E-04 & 	1.313 &	4.22E-03 &	  	1.205   \\
        & 128	& 7.13E-04 & 1.928 &1.79E-04 & 1.803	 &	1.07E-03 &			1.978  \\
        & 256	& 1.80E-04 & 1.986 &4.64E-05 & 1.950	 &	2.84E-04 &			1.913  \\
        \hline
        \multirow{4}{*}{$10^{-2}$} & 32	& 1.71E-03	& ---    & 1.96E-04 &  ---   & 	9.98E-02 & ---        \\
        & 64	& 5.37E-04	&1.668& 6.65E-05 &1.558& 	3.69E-02 &	1.436 \\
        & 128	& 1.42E-04	&1.920& 1.83E-05 &1.862& 	8.78E-03 &	2.069 \\
        & 256	& 3.54E-05	&2.003& 4.69E-06 &1.963& 	2.33E-03 &	1.915 \\
        \hline
        \multirow{5}{*}{$10^{-3}$} & 32	&   1.48E-03	&   ---    & 	2.00E-05	& ---    & 1.04E+00 &	---        \\
        & 64	&   4.11E-04	& 1.845 & 	6.75E-06	&1.566& 3.70E-01 &			1.485\\
        & 128 & 	1.07E-04	& 1.943 & 	1.83E-06	&1.881& 8.73E-02 &			2.083\\
        & 256 & 	2.63E-05	& 2.023 & 	4.69E-07	&1.964& 2.31E-02 &			1.918\\
        \hline
    \end{tabular}
    \caption{Travelling isentropic vortex: $L^1$ error and convergence rates at $t_f = 0.2$ for $M=1$ and $M=10^{-2}$ under moderate CFL numbers on Cartesian grid composed of $N\times N$ cells. }
    \label{tab:EOCEuler}
\end{table}

\begin{table}[htpb]
    \centering
    \renewcommand{\arraystretch}{1.25}
    \begin{tabular}{cccccc}
        $M$ & $10^{0}$ & $10^{-1}$ & $10^{-2}$& $10^{-3}$ & $10^{-4}$\\ \hline
        $\rho$ &  5.26E-04 &
        6.49E-04 &
        4.19E-04 &
        9.98E-06 &
        9.21E-08
        \\
        $p$ & 3.83E-04 &
        3.66E-04 &
        3.51E-05 &
        6.44E-07 &
        8.05E-09
        \\\hline
    \end{tabular}
    \caption{Isentropic vortex: relative $L^1$ errors of density and pressure for different Mach number regimes on a fixed grid with $N=64$ with respect to $\rho_0, p_0$ and $\nu_{mat} = 0.1$ in the CFL condition at final time $t_f = 0.1$. }
    \label{tab:EulerAP}
\end{table}

\subsubsection{Riemann Problems}
Next, we consider three classical Riemann Problems (RPs). The first one is the Sod shock tube problem (RP1) which is situated in the compressible flow regime with $M\approx 0.9$ on the contact wave.
Then we consider a slow moving contact problem (RP2) with $M\approx 6\cdot 10^{-3}$ on the contact wave and finally a contact wave moving with constant velocity under constant pressure while the density changes 5 orders of magnitude (RP3).
The set-up for all three RPs is summarized in Table \ref{tab:RPEuler} which is taken from \cite{Sod1978, AbbIolPup2017, Boscheri2021}.
The CFL condition for all test cases is given by \eqref{eq:CFL} with the material CFL number $\nu = 0.5$

We compare our results obtained with the new \IMEX scheme against the exact Riemann solution from \cite{Toro2009}.
For for RP1 and RP2, we also compare with the fully IMplicit Finite Volume Predictor Corrector scheme \IM which is derived from the relaxation model \eqref{sys:Relax} where all flux terms are discretized implicitly.
For details, see \cite{Thomann2023}.

\begin{figure}[htpb]
    \includegraphics[width=0.34\textwidth]{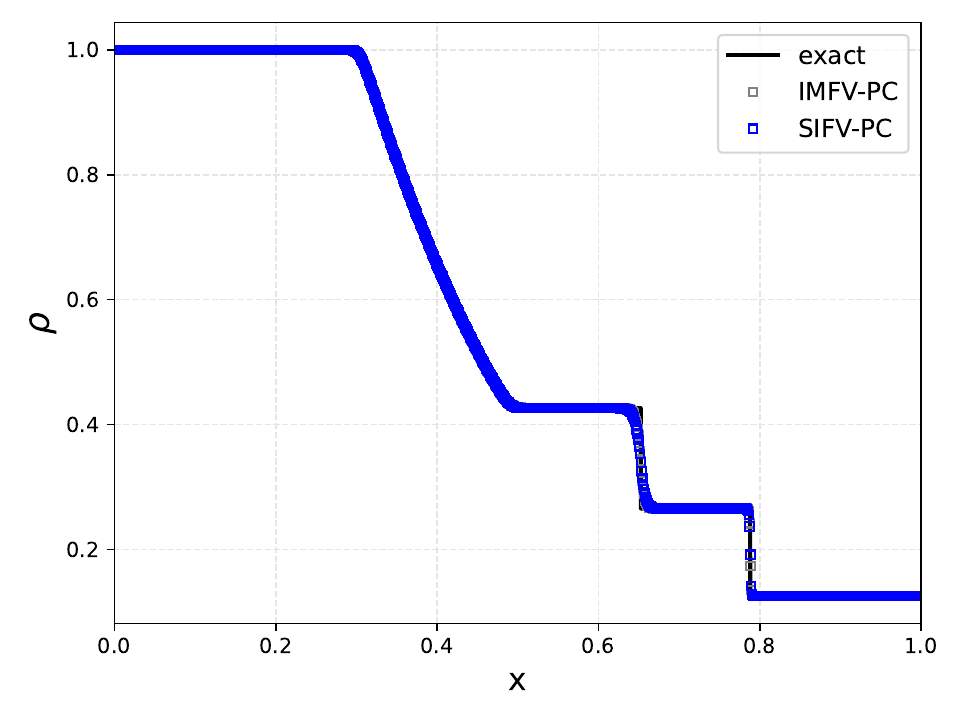}
    \includegraphics[width=0.34\textwidth]{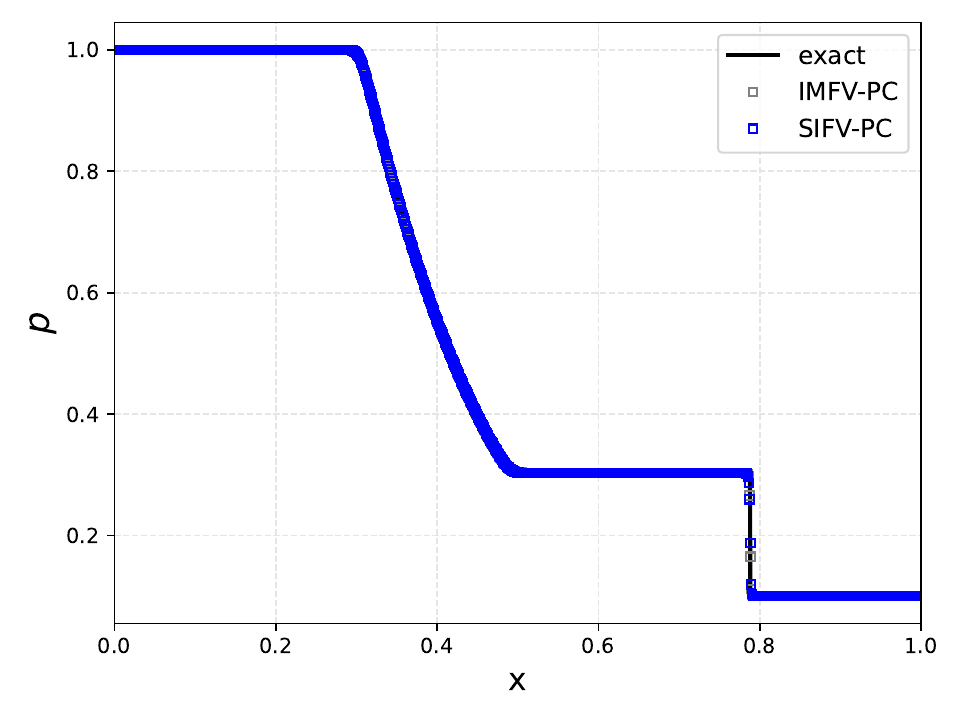}
    \includegraphics[width=0.34\textwidth]{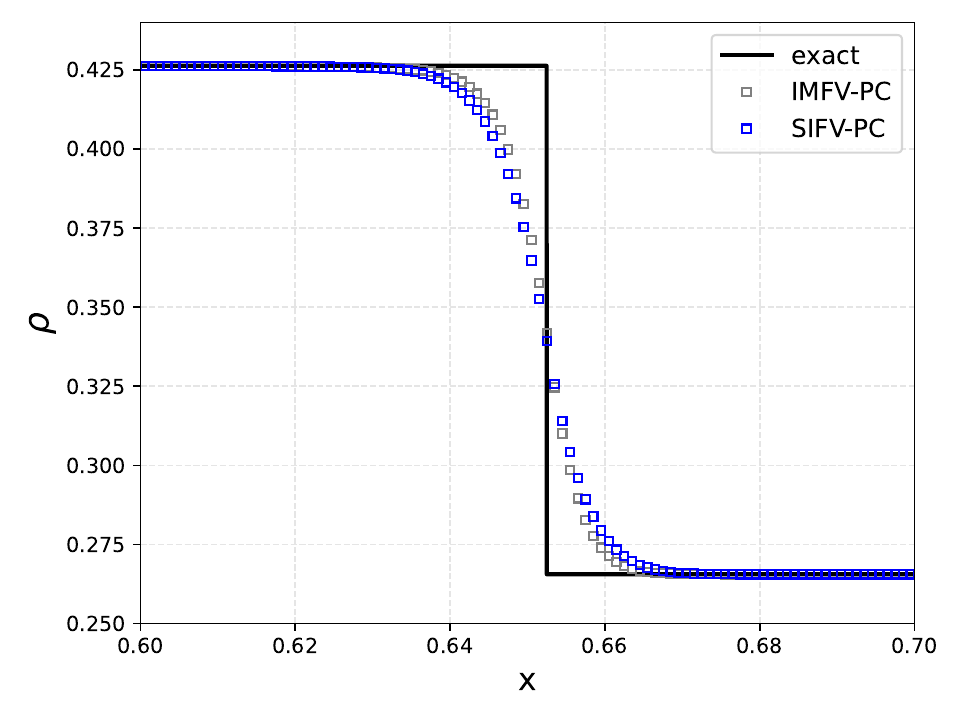}
    \includegraphics[width=0.34\textwidth]{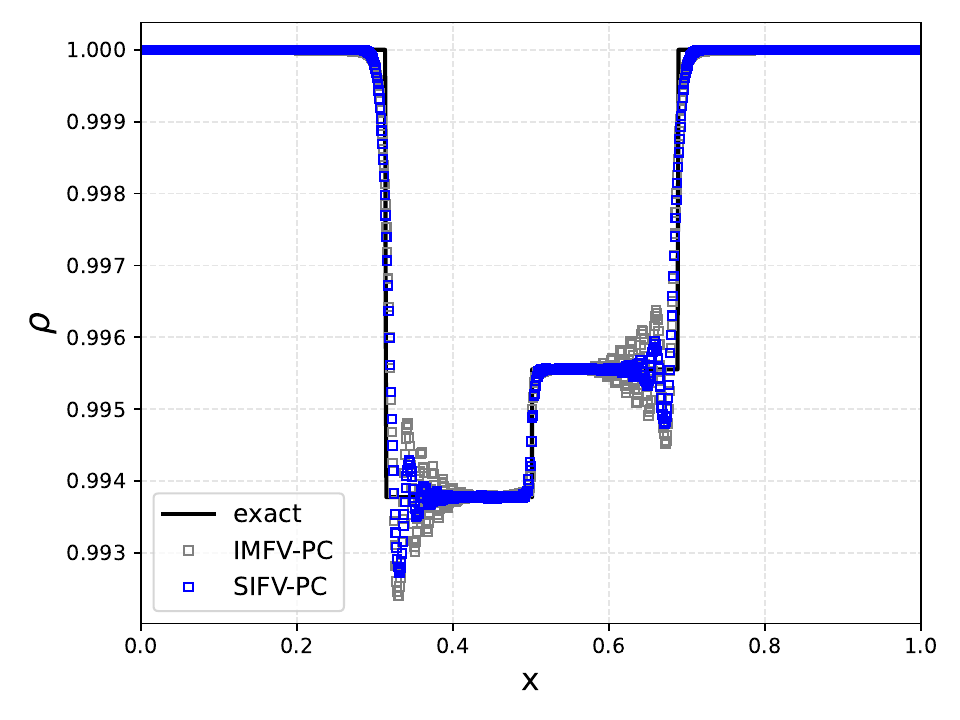}
    \includegraphics[width=0.34\textwidth]{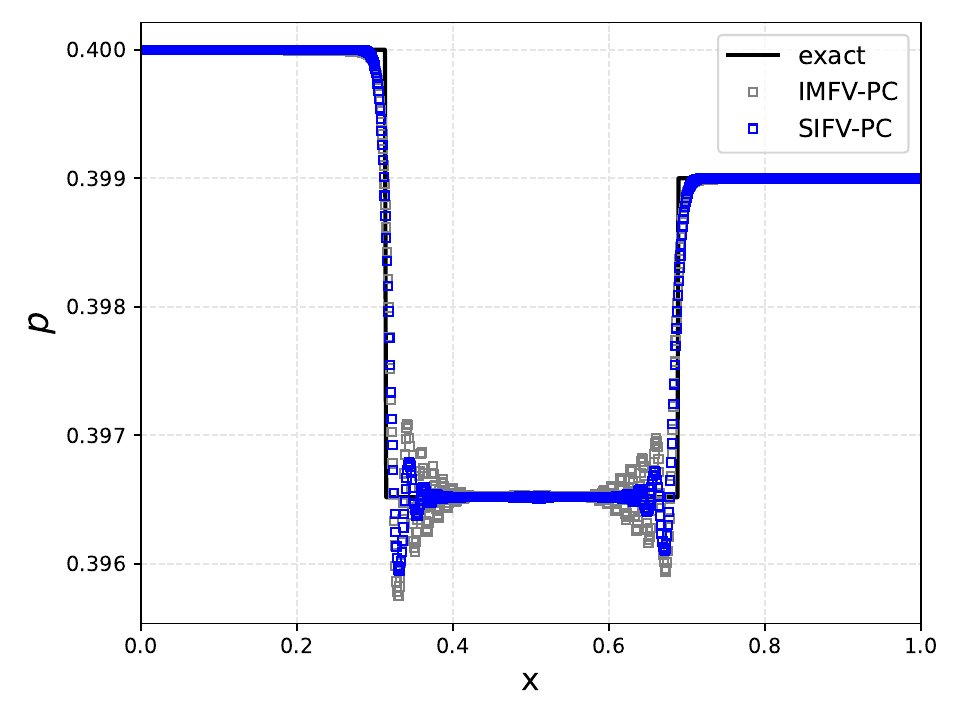}
    \includegraphics[width=0.34\textwidth]{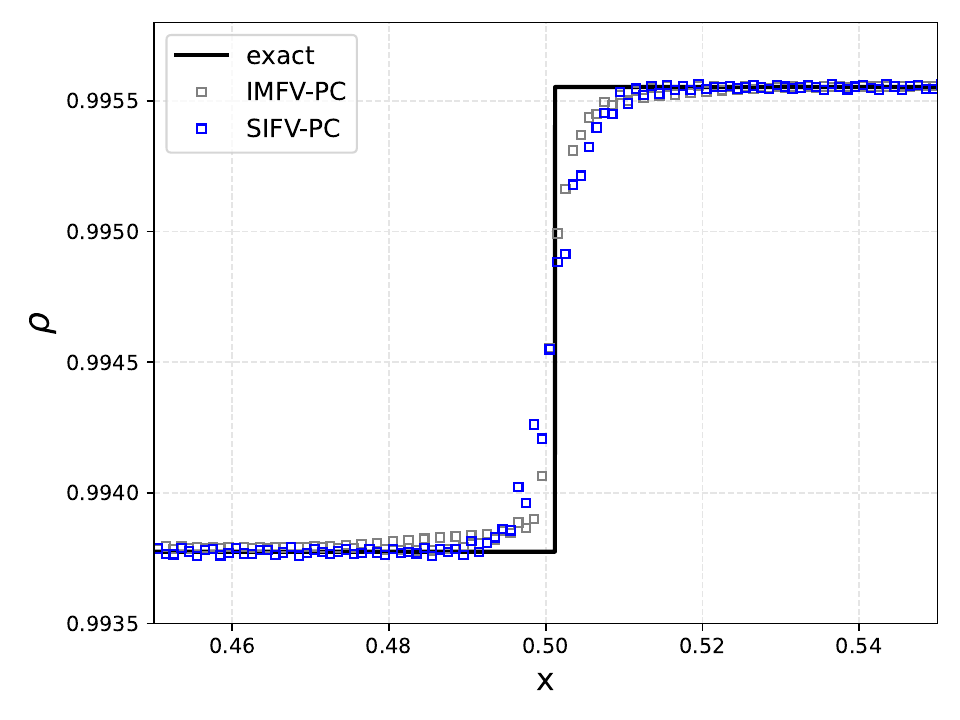}
    \caption{Euler Riemann Problems RP1 and RP2 from Table \ref{tab:RPEuler}: Sod and slow moving contact test. Comparison of the new second order semi-implicit \IMEX scheme against the fully implicit \IM scheme from \cite{Thomann2023} and the exact solution on 1000 cells. }
    \label{fig:RP12Euler}
\end{figure}

\begin{figure}[htpb]
    \includegraphics[width=0.34\textwidth]{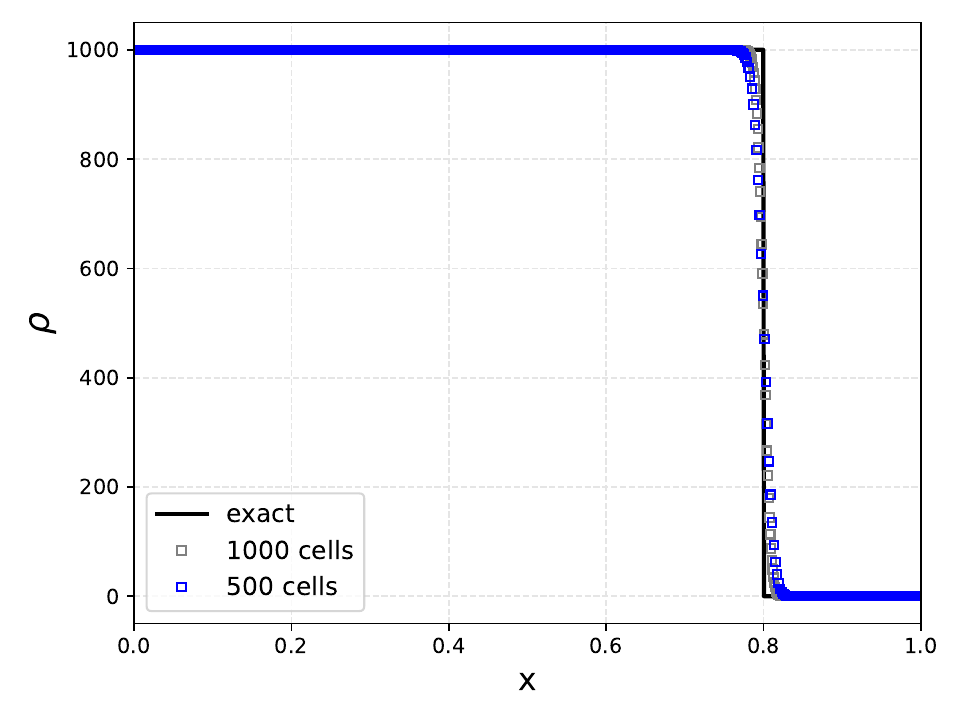}
    \includegraphics[width=0.34\textwidth]{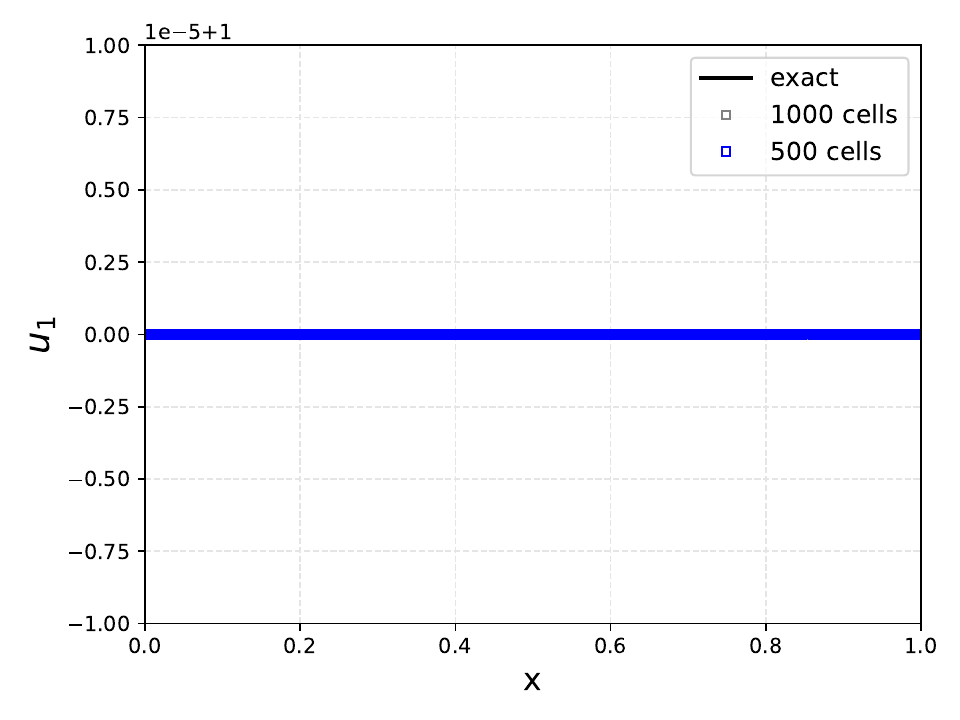}
    \includegraphics[width=0.34\textwidth]{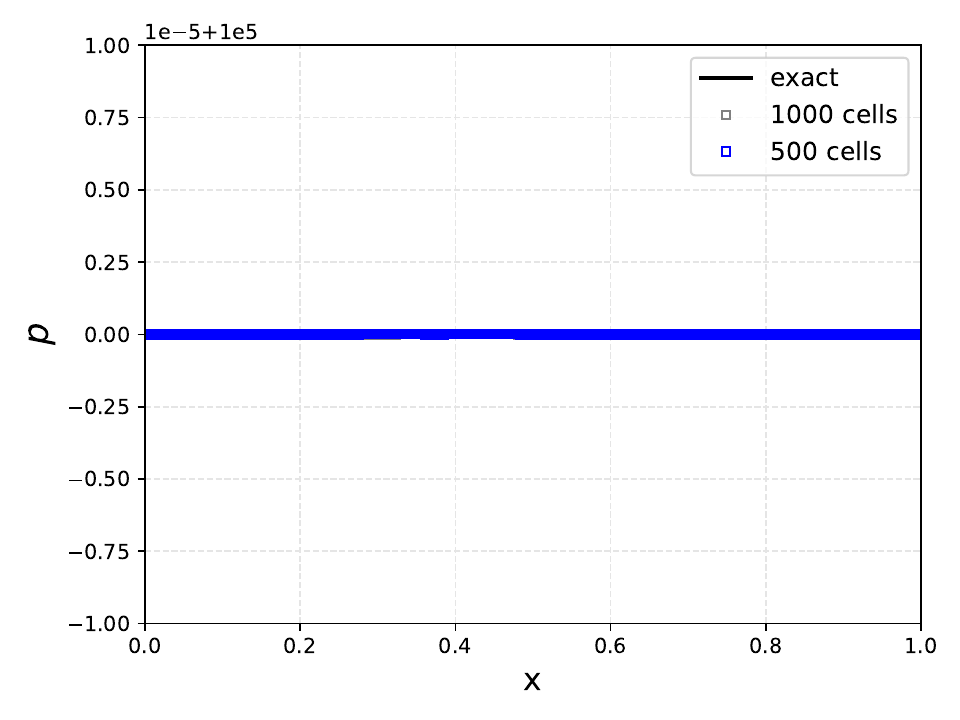}
    \caption{Euler Riemann Problem RP3 from Table \ref{tab:RPEuler}: Travelling contact wave at $t_f= 0.5$ computed with the first order \IMEX scheme on two grids with $N=500$ and $N=1000$ cells. }
    \label{fig:RP3Euler}
\end{figure}

The numerical solutions are depicted in Figures \ref{fig:RP12Euler} and \ref{fig:RP3Euler} and verify the ability of the new scheme to accurately capture the contact waves in different flow regimes.
However, in RP2 we observe oscillations close to the left and right travelling shock waves.
This is due to an under-resolution of the fast travelling waves in time due to the use of a material time step combined with the low diffusivity of the numerical scheme as the maximal velocity is $8\cdot 10^{-3}$.
The oscillations can be avoided or reduced by choosing a smaller time step and increasing the numerical diffusion.
This however leads to a more diffusive contact wave, thus a trade-off between contact and shock-wave resolution has to be made.

In Figure \ref{fig:RP3Euler}, the results for RP3 for the first order \IMEX scheme are depicted.
It verifies the contact-preservation property proven in Lemma \ref{lem:contact} proven for the first order scheme which is obtained due to the Toro-Vazques splitting \cite{Toro2012}.
The fully implicit \IM scheme does not have this property and is thus not included in the plot.

\begin{table}[htpb]
    \centering
    \renewcommand{\arraystretch}{1.25}
    \begin{tabular}{l c c c c cccc}
        RP & $\rho_L$ &
         $\rho_R$ & $u_L$ & $u_R$ & $p_L$ & $p_R$ & $x_0$ & $t_f$ \\\hline
       RP1 & 1.0 & 0.125 & 0.0 & 0.0 & 1.0 & 0.1 & 0.5 & 0.1644 \\
       RP2 & 1.0 & 1.0 & 0.0 & 0.008 & 0.4 & 0.399 & 0.5 & 0.25 \\
       RP3 & $10^3$ & 0.01 & 1.0 & 1.0& $10^5$ & $10^5$ & 0.3 & 0.5 \\\hline
    \end{tabular}
\caption{Initialization of Riemann Problems for the Euler equations. The initial conditions are given in terms of left (L) and right (R) states separated by an initial discontinuity at $x = x_0$ on the computational domain $[0,1]$. Due to the one-dimensional set-up, for the velocities in $y$- and $z$- direction we set $u = u_x$ and $u_y,u_z = 0$. Further, the final time is denoted by $t_f$ and in all cases $\gamma = 1.4$. }
\label{tab:RPEuler}
\end{table}

\subsubsection{Kelvin-Helmholtz}
As final test case for the Euler equations, we simulate a Kelvin-Helmholtz instability in the set-up of \cite{Feireisl2021}.
The initial condition for the KH problem on the unit box $~\Omega = [0,1]\times [0,1]$ is given by
\begin{equation}
    \label{eq:KH}
    \left(\rho, u_1, u_2, p\right) =
    \begin{cases}
        \left(2, -0.5, 0, 2.5\right) &\text{if } I_1 < x_2 < I_2 \\
        \left(1,  0.5, 0, 2.5\right) &\text{else.}
    \end{cases}
\end{equation}
The interface profiles, denoted by $I_j = I_j(\xx)$ for $j=1,2$, are defined as $I_j = J_j + \varepsilon Y_j(\xx)$ which generate small perturbations around the constant horizontal lines $J_1 = 0.25$ and $J_2 = 0.75$.
Further, we set
\begin{equation}    Y_j(\xx) = \sum_{k=1}^m a_j^k \cos(b_j^k + 2 k \pi x_1), \quad j= 1,2,
\end{equation}
where the parameters $a_j^k \in [0,1]$ and $b_j^k$, $j=1,2$, $k=1, \dots, m$ are fixed.
The coefficients $a_j^k$ are normalized such that $\sum_{k=1}^m a_j^k = 1$ to guarantee that the perturbation around $J_j$ are of the order of $\varepsilon$, thus $|I_j(\xx) - J_j| \leq \varepsilon$ for $j=1,2$.
In the numerical simulation we use $m=10$ modes and the perturbation $\varepsilon = 0.01$, as done in \cite{Feireisl2021}.
Further we deploy periodic boundary conditions in both directions and set $\nu = 0.25$ in the material CFL condition \eqref{eq:CFL}.

\begin{figure}[htpb]
    \includegraphics[width=0.36\textwidth]{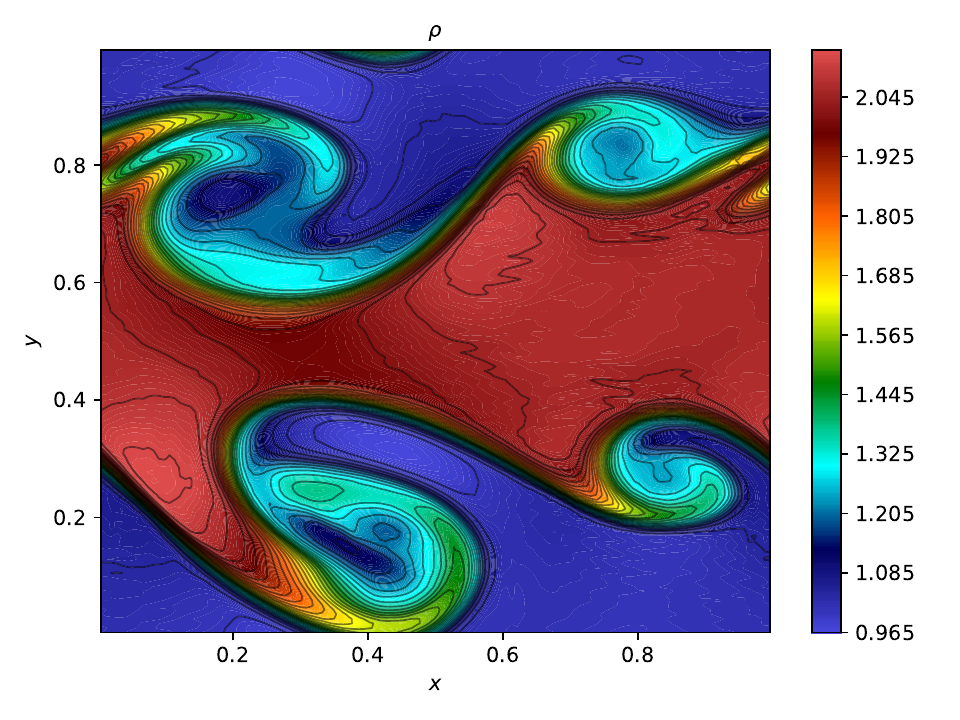}\hskip-2mm
    \includegraphics[width=0.36\textwidth]{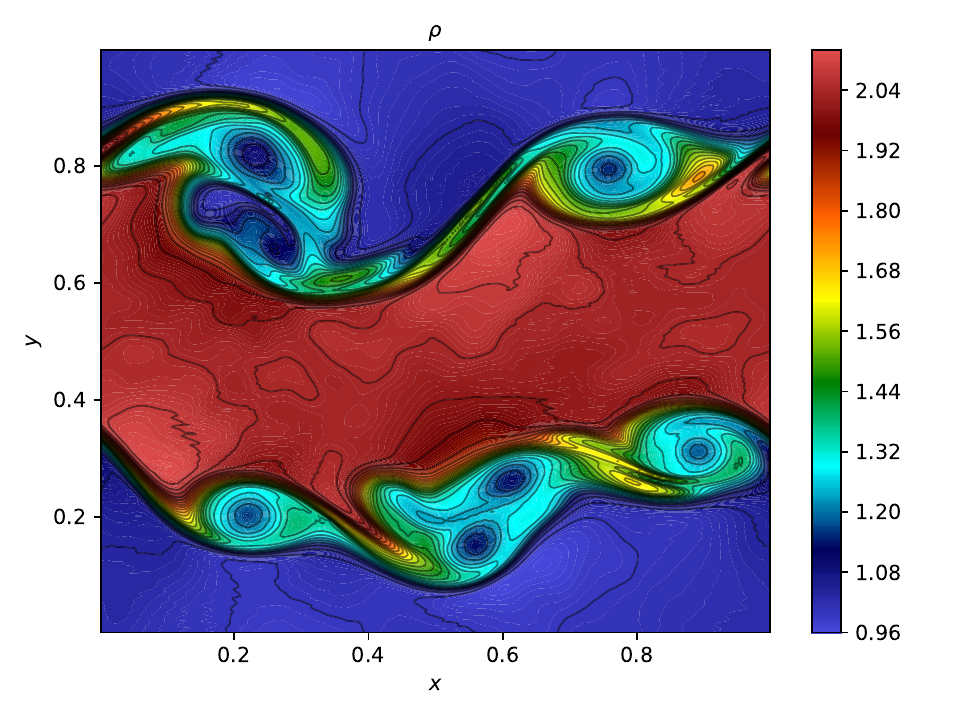}\hskip-3mm
    \includegraphics[width=0.36\textwidth]{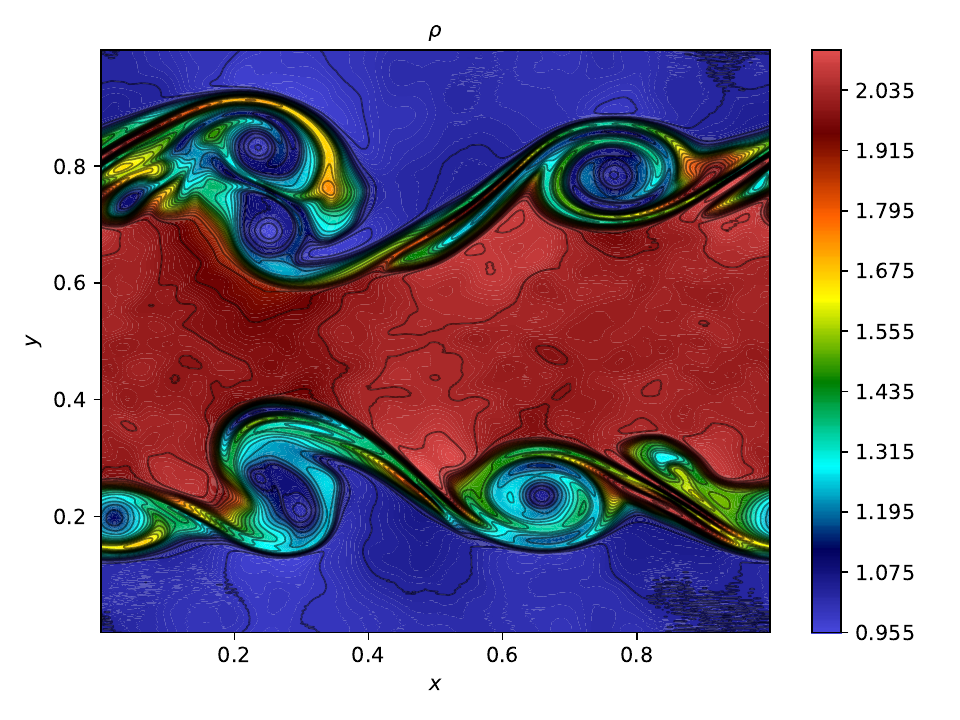}
    \caption{Kelvin-Helmholtz instability: Density $\rho$ obtained with the second order \IMEX scheme computed with $N\times N$ cells in the Mach number regime $M\approx 0.25$ at $t=2$ displayed with 31 equidistant contour lines. Left: $N=128$. Center: $N=256$. Right: $N=512$.}
    \label{fig:KH}
\end{figure}

In Figure \ref{fig:KH} the density obtained with three resolutions $N=128,256,512$ is displayed.
We see that as the mesh is refined, the formed vortices become more defined.
The qualitative difference of the solutions is due to the initial perturbations of the interface profiles.

\subsection{The ideal MHD equations}
Next we consider the ideal MHD equations in the splitting \eqref{eq:MHDsplit} using the three-split scheme.
To complete the description of the numerical scheme, we have to define the matrices $\AA_1$ and $\AA_2$ which for all test cases are given by $\AA_1 = a_1 \mathbb{I}$ with $a_1 = \max_i\max_j|\lambda_j^{1,\fast}(\qq_i^n)|$ and $\AA_2 = \max_i\max_j|\lambda_j^{2,\fast}(\qq_i^n)|$.

\subsubsection{Travelling Balsara vortex}

We test the accuracy of the numerical scheme on the analytical vortex solution forwarded in \cite{Balsara2004} in the version of \cite{Boscheri2024} on the two-dimensional computational domain for $\xx = (x_1, x_2)^T \in \Omega=[-5;5]\times[-5;5]$.
The vortex is constructed as a perturbation $\delta \qq$ on top of an unperturbed magnetohydrodynamics background flow given by the constant states $\qq_0 = (\rho_0, \uu_0, p_0, \mathbf{0}, 0)$.
The perturbations of the velocity and magnetic field are given by
\begin{eqnarray}
    &&(\delta u_1, \delta u_2) = \frac{\tilde{v}}{2 \pi} \, \exp\left( \frac{1- r^2}{2}\right) \cdot (- r \sin(\theta), r \cos(\theta)), \\
    &&(\delta B_1, \delta B_2) = \frac{\tilde{B}}{2 \pi}  \, \exp\left( \frac{1- r^2}{2}\right) \cdot (- r \sin(\theta), r \cos(\theta)),
\end{eqnarray}
with the constants $\tilde{u} = \sqrt{2\pi}$ and $\tilde B = \sqrt{4 \pi}$.
The perturbation of the magnetic field is generated by the magnetic potential perturbation in $z$-direction $\delta A_z = \frac{\tilde B}{2\pi} \, \exp(1-r^2)$.
The pressure is set such that it balances the centrifugal force generated by the circular motion against the centripetal force stemming from the tension in the magnetic field lines. In radial coordinates, with the generic radial position $r = \sqrt{x_1^2 + x_2^2}$, the equation for the pressure perturbation writes
\begin{equation}
    \label{eq:ODE_vortex_dp}
    \frac{\partial \delta p}{\partial r} = \left( \frac{\rho u_\theta^2}{r} - \frac{\delta B_\theta^2}{4 \pi r} - \frac{\partial (\delta B_\theta)^2}{8 \pi r}\right),
\end{equation}
with the definitions
\begin{equation}
    \delta u_\theta = \frac{\tilde{u}}{2 \pi} \, r \, \exp\left( \frac{1- r^2}{2}\right), \qquad \delta B_\theta = \frac{\tilde{B}}{2 \pi} \, r \, \exp\left( \frac{1- r^2}{2}\right).
\end{equation}
Integration of \eqref{eq:ODE_vortex_dp} yields
\begin{equation}
    \label{eq:Vortex_dp}
    \delta p = \frac{1}{2} e^{1-r^2} \left(\frac{1}{8 \pi}\frac{\tilde{B}^2}{(2 \pi)^2} \left(1-r^2\right)-\frac{1}{2}\rho \left(\frac{\tilde{u}}{2\pi}\right)^2\right).
\end{equation}
We set the background velocity field to $\uu_0 = (1,1)^T$, thus the vortex is moving diagonally through the computational domain, and the background pressure is $p_0 = 1$.
To vary the Mach and Alfv\'en number regime, we lower the density as suggested in \cite{Boscheri2024} to obtain faster acoustic and Alfv\'en speeds.
To that end we consider two different background densities, $\rho_0 = 1$ and $\rho = 10^{-2}$.
The simulation are run on a sequence of successively refined meshes until the final time $t_f = 0.2$ such that the solution does not leave the prescribed computational domain.
The material CFL number $\nu = 0.25$ in the CFL condition \eqref{eq:CFL} is applied.
In Table \ref{tab:EOCMHD} the $L^1$ error in $\rho, |\uu|, p, |\BB|$ are given.
In addition, we give the maximum sound speed $c$ and Alfv\'en speed $b$.
We observe for the considered regimes, the expected EOC of order two is achieved.

\begin{table}[htbp]
    \centering
    \renewcommand{\arraystretch}{1.25}
    \begin{tabular}{cccrcccccccc}
        $\rho_0$ & $c$ & $b$ & $N$ & $|\uu|$ & & $\rho$ & & $|\BB|$ & & $p$ & \\ \hline
        \multirow{4}{*}{$10^0$} & \multirow{4}{*}{1.18} &\multirow{4}{*}{$1.59\cdot 10^{-1}$}&32&3.98E-02 & ---   &	9.66E-03 &	---   &	7.80E-02 & ---    &	2.97E-02 &	  ---  \\
        &&&64&7.34E-03 & 2.439 &	2.80E-03 &	1.786	&	1.76E-02 & 2.151	&	6.76E-03 &	2.135  \\
        &&&128&1.58E-03 & 2.218 &	7.17E-04 &	1.966	&	4.24E-03 & 2.048	&	1.64E-03 &	2.039  \\
        &&&256&3.84E-04 & 2.036 &	1.80E-04 &	1.997	&	1.07E-03 & 1.989	&	4.17E-04 &	1.978  \\

\hline
        \multirow{4}{*}{$10^{-2}$} & \multirow{4}{*}{$1.20\cdot 10^1$} &\multirow{4}{*}{1.59}&32	  &8.80E-02 & ---  	& 9.38E-05 & ---   & 	5.14E-02 & ---   & 	6.78E-03	& ---   \\
        &&&64	  &5.70E-02 & 0.626	& 7.75E-05 & 0.275 & 	3.40E-02 & 0.595 & 	2.78E-03	& 1.285 \\
        &&&128	&1.45E-02 & 1.977	& 2.10E-05 & 1.882 & 	1.20E-02 & 1.499 & 	9.51E-04	& 1.548 \\
        &&&256	&3.44E-03 & 2.073	& 5.13E-06 & 2.033 & 	3.00E-03 & 2.006 & 	2.11E-04	& 2.174 \\\hline
    \end{tabular}
\caption{MHD vortex: $L^1$ errors and convergence rates at $t_f = 0.2$ for the background densities $\rho = 1, 10^{-2}$ under varying sound speeds $c$ and Alfv\'en speed $b$ on Cartesian grid composed of $N\times N$ cells.  }
\label{tab:EOCMHD}
\end{table}

\subsubsection{Riemann Problems}
\begin{figure}[htpb]
        \includegraphics[width=0.34\textwidth]{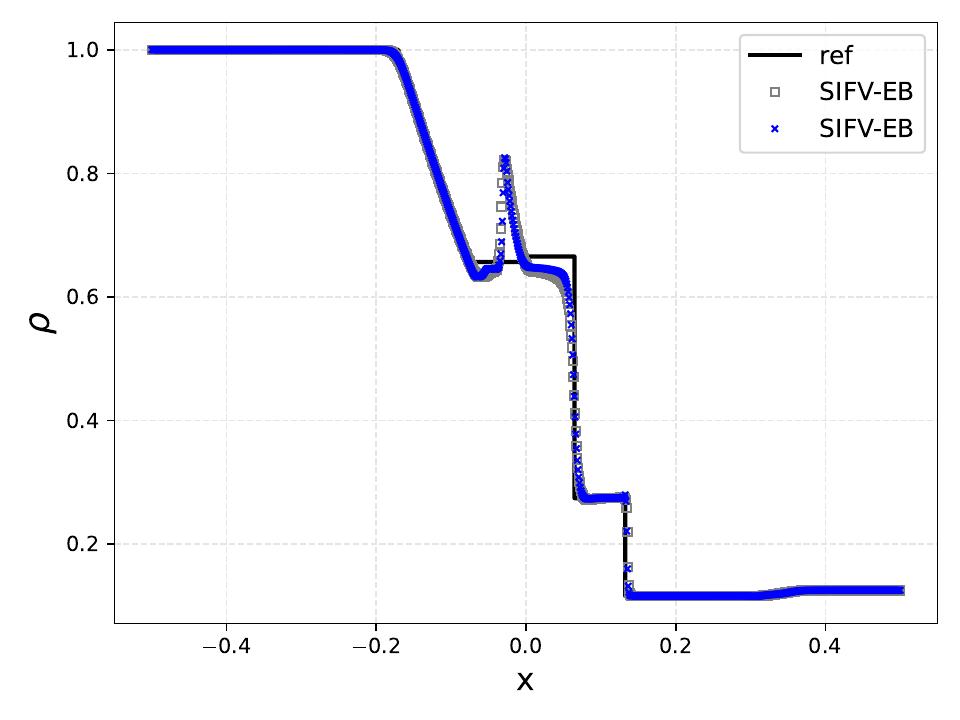}
        \includegraphics[width=0.34\textwidth]{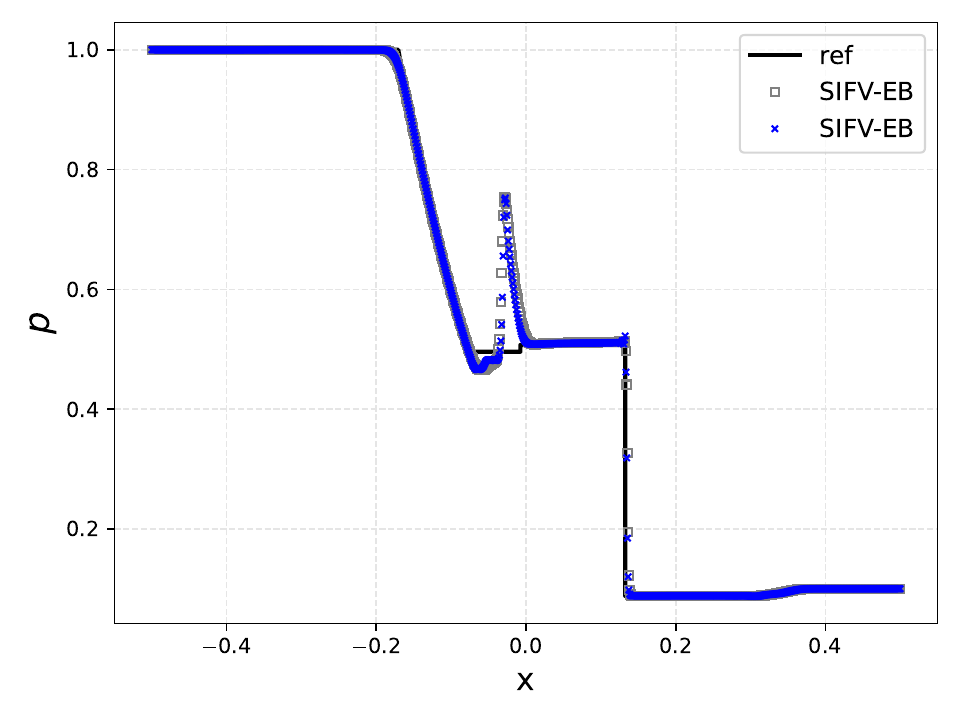}
        \includegraphics[width=0.34\textwidth]{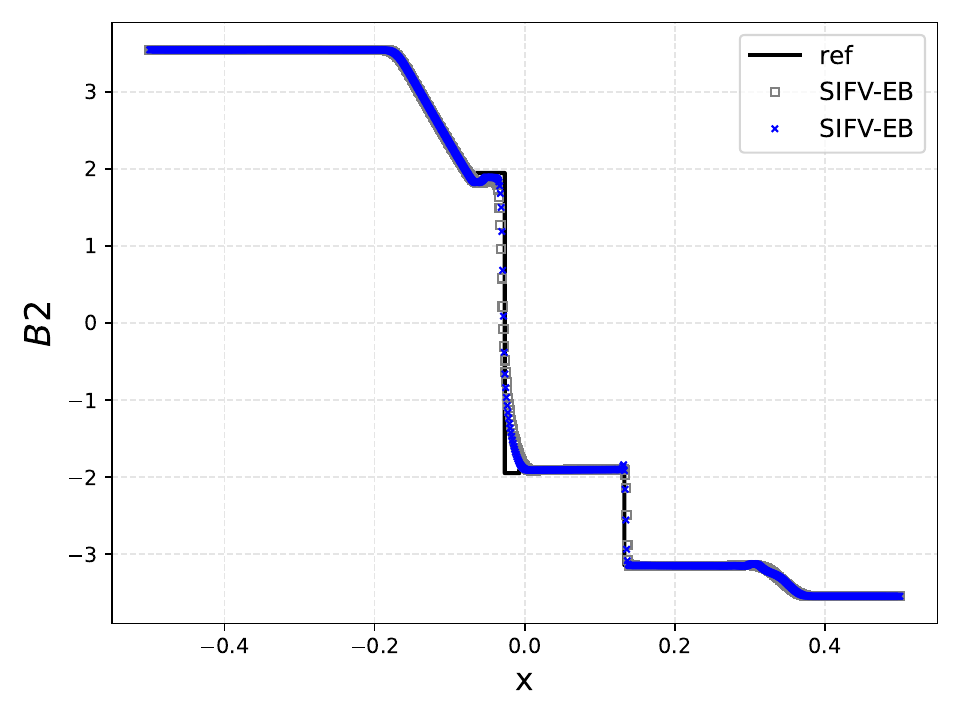}

        \includegraphics[width=0.34\textwidth]{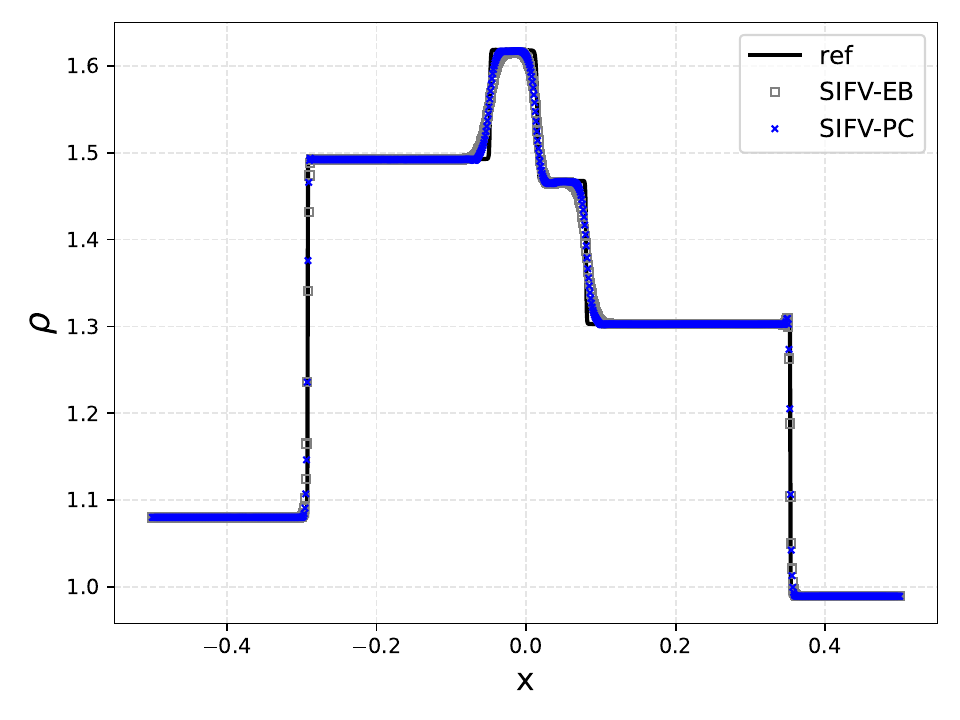}
        \includegraphics[width=0.34\textwidth]{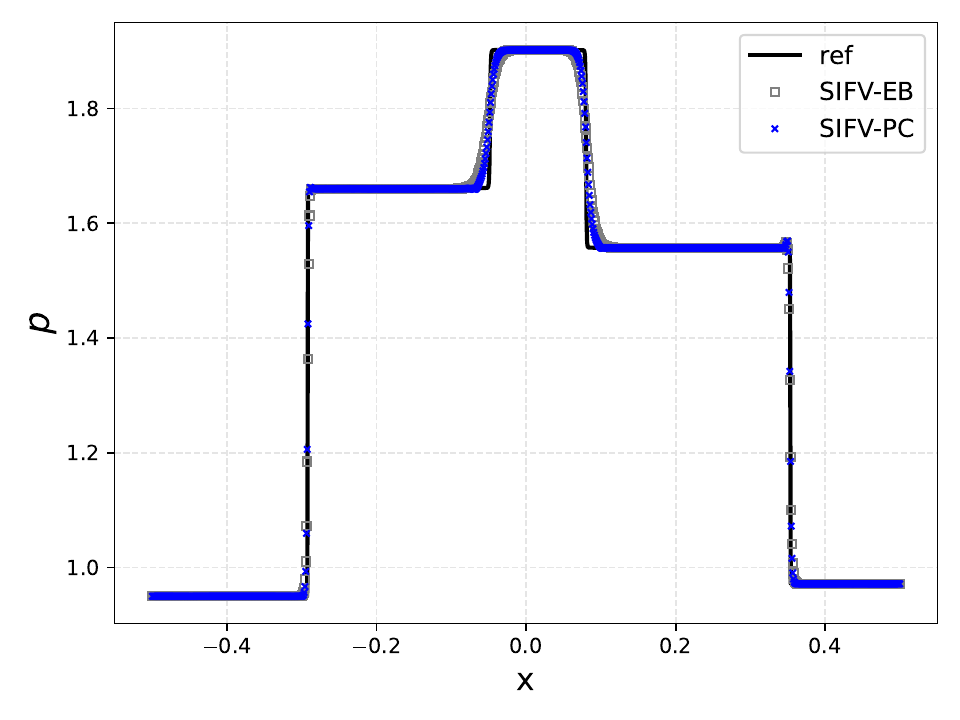}
        \includegraphics[width=0.34\textwidth]{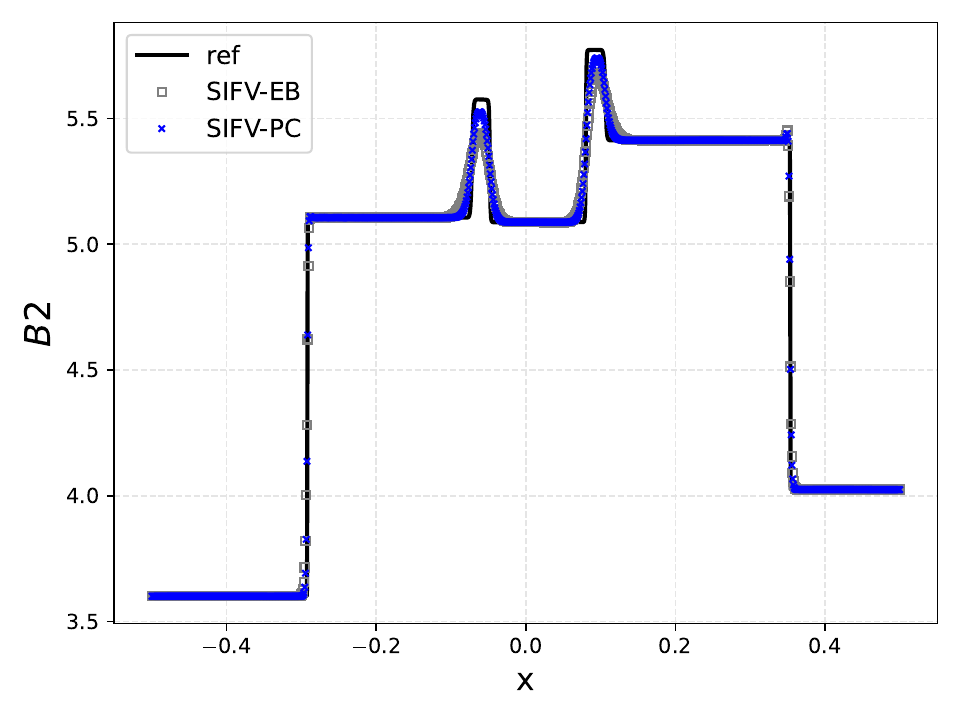}

        \includegraphics[width=0.34\textwidth]{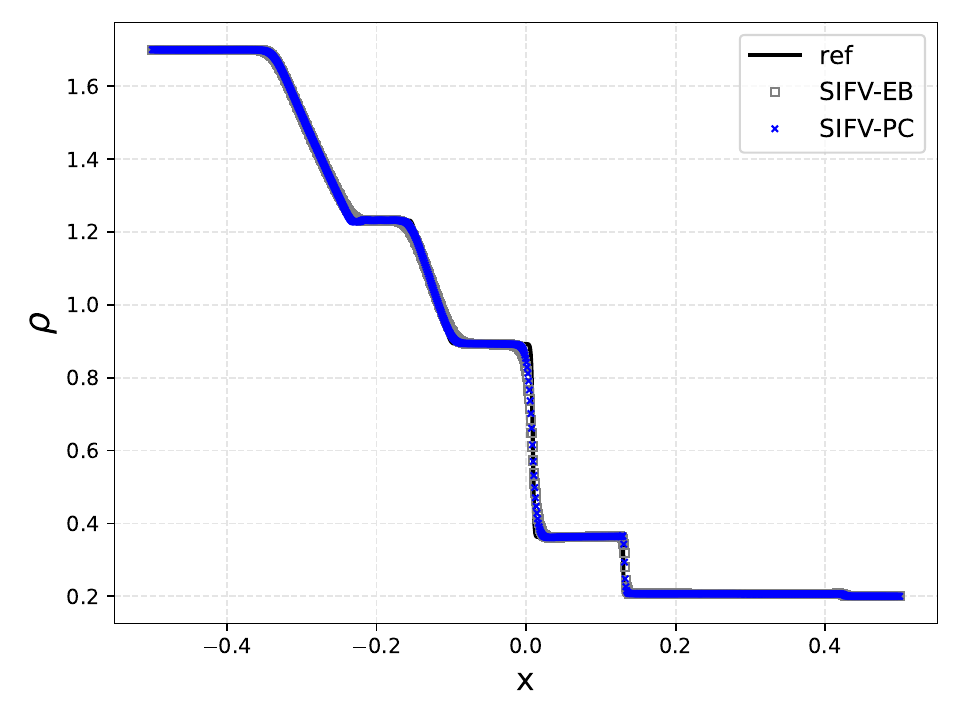}
        \includegraphics[width=0.34\textwidth]{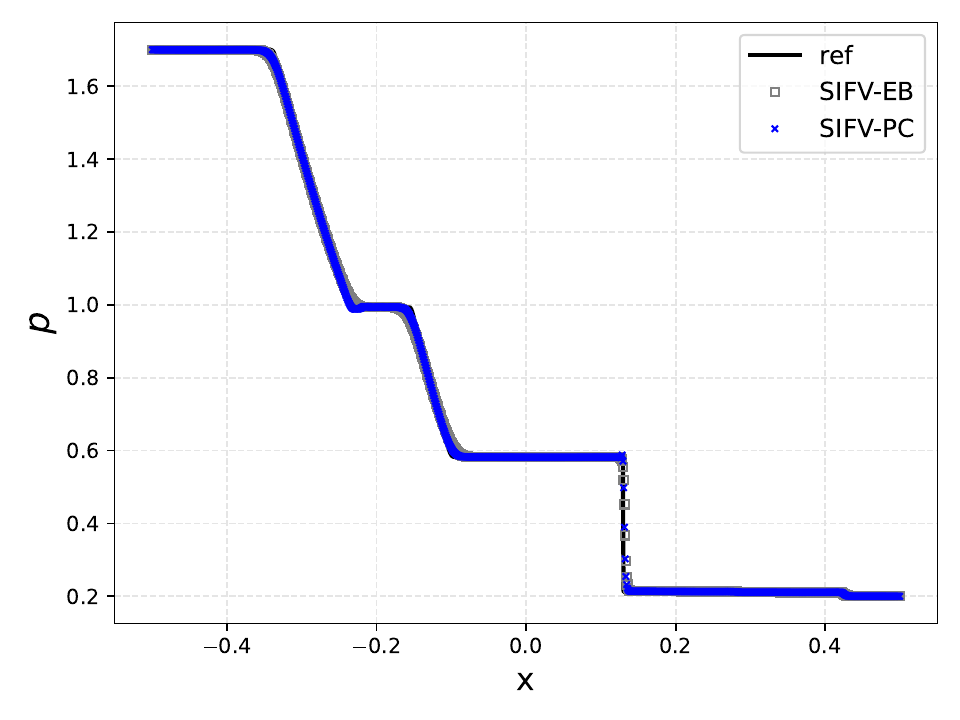}
        \includegraphics[width=0.34\textwidth]{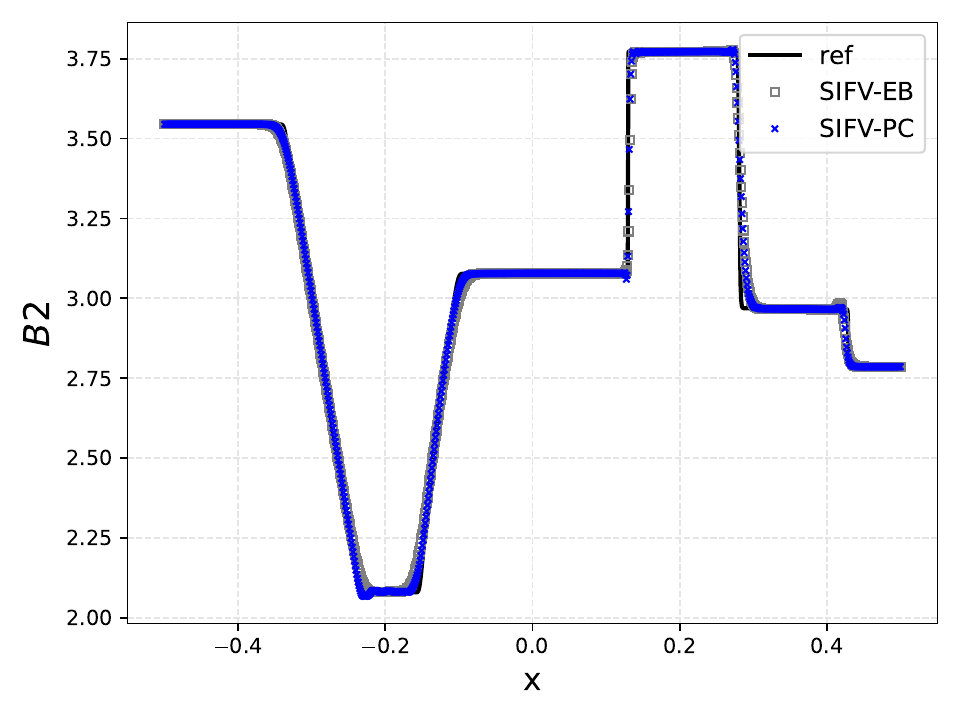}

        \includegraphics[width=0.34\textwidth]{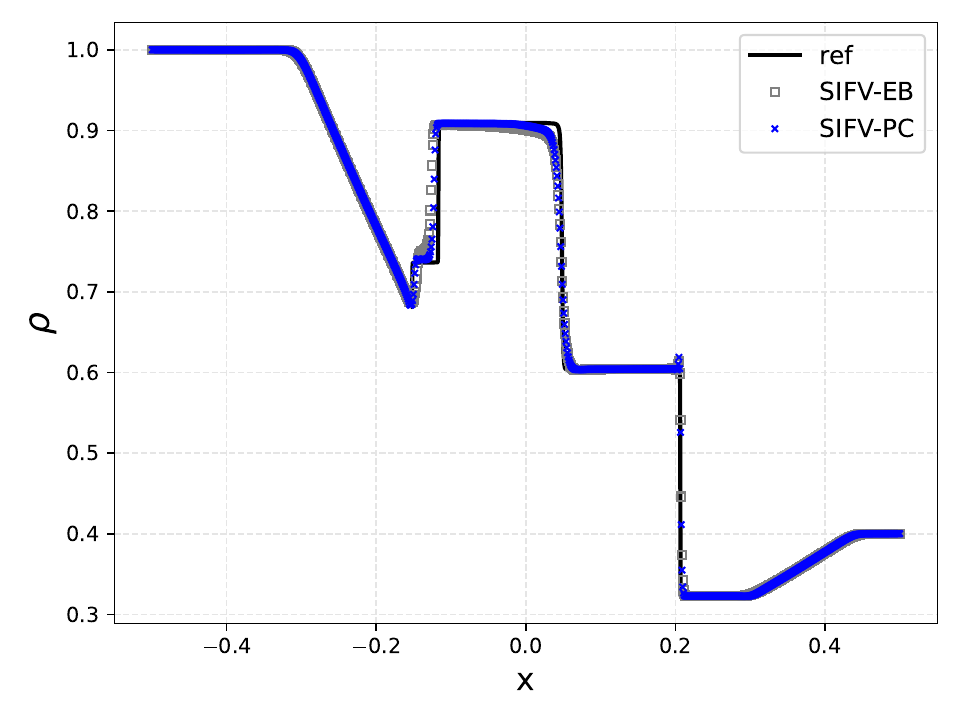}
        \includegraphics[width=0.34\textwidth]{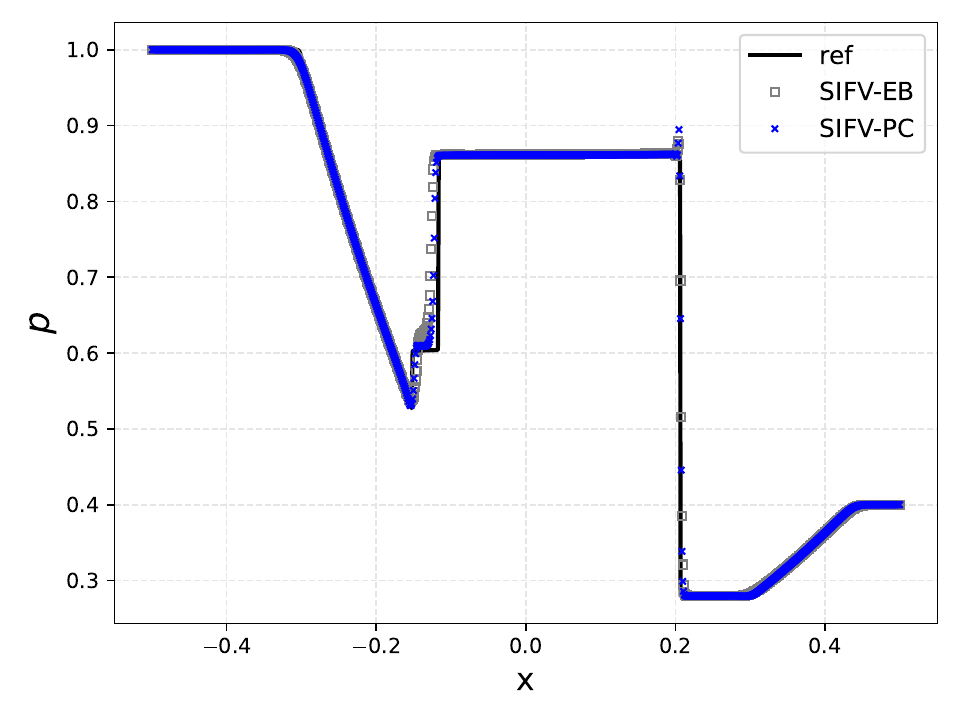}
        \includegraphics[width=0.34\textwidth]{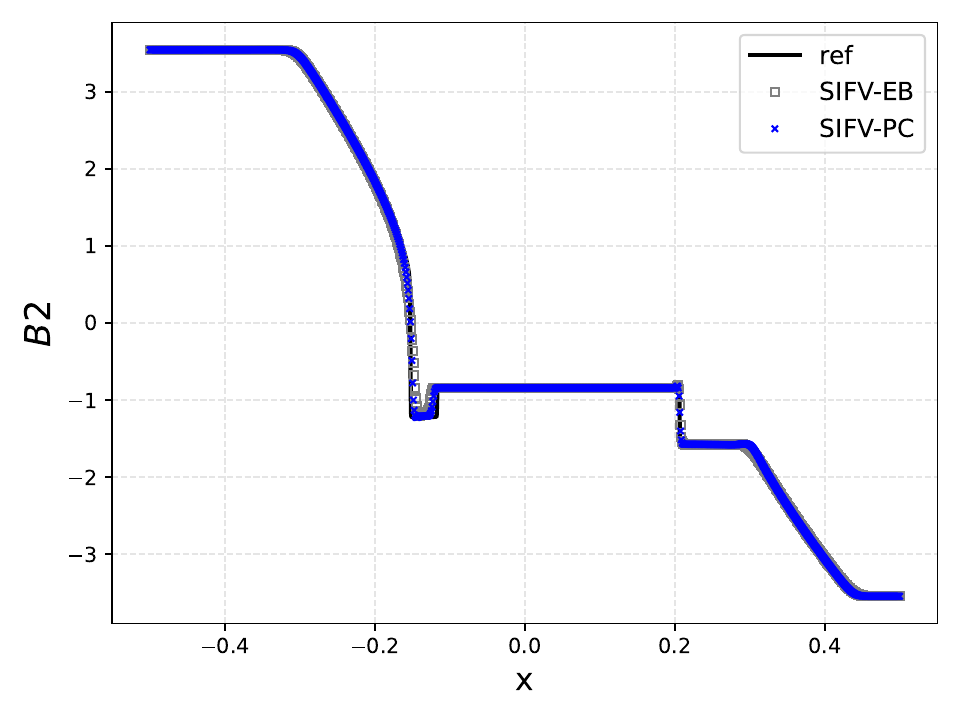}

        \caption{MHD Riemann Problems RP1-4 from Table \ref{tab:RPMHD}. Comparison with the second-order semi-implicit SIFV-EB scheme from \cite{Boscheri2024} on 1000 cells.}
\label{fig:RPMHD}
\end{figure}

Next, we consider a series of RPs for the one-dimensional ideal MHD equations.
They are standard test cases and the set-up is summarized in Table \ref{tab:RPMHD} taken from \cite{Dumbser2016}.
For instance, they have been used to validate for instance the Osher-type scheme \cite{Dumbser2011} or the HLLEM Riemann solver presented in \cite{Dumbser2016}.
We compare the results obtained with the new \IMEX scheme against the exact solution, that has been computed with an exact Riemann solver that has kindly been provided by Falle and Komissarov \cite{Falle2002,Falle2001} and which is labelled in Figure \ref{fig:RPMHD} as reference solution \textit{ref}.
For an alternative exact Riemann solver of the MHD equations, see also \cite{torrilhon2003}.
Further we compare to a semi-implicit three-split scheme for the MHD equations recently forwarded in \cite{Boscheri2024} denoted by the authors as \SIMHDthree scheme.
It employs a similar strategy than the three-split schemes introduced in this work, but the linearity of the \SIMHDthree scheme relies on the linearity of the equation of state with respect to the total energy and moreover is based on a different splitting than the Fambri-splitting \cite{Fambri2021} employed in this work.
These test cases verify the robustness of the new \IMEX scheme in the case of moderate and high Mach number flows. The computational domain is given by $\Omega=[0,1]$ and the grid is made of 1000 cells along the $x-$direction with Neumann boundary conditions.
The material CFL number in \eqref{eq:CFL} is set to $\nu = 0.5$.
The novel \IMEX scheme behaves similar to the \SIMHDthree scheme as it has the same time step restriction and similar diffusivity and is able to capture all the waves in the considered RPs.
In particular, the overshoot in density and pressure in RP1 is an artifact that also appears in approximate solutions produced with the HLLEM scheme from \cite{Dumbser2016}, the Osher-type scheme \cite{Dumbser2011} or the Rusanov scheme.

\begin{table}[htpb]
    \centering
    \renewcommand{\arraystretch}{1.25}
    \begin{tabular}{l c c c cccccccc}
        RP & &$\rho$ & $u_x$ & $u_y$ & $u_z$ & $p$ & $B_x$ & $B_y$ & $B_z$ & $x_0$ & $t_f$ \\\hline
        \multirow{2}{*}{RP1} & L: & 1.0 & 0.0 & 0.0 & 0.0 & 1.0 & 0.75$\sqrt{\mu_0}$ & $\sqrt{\mu_0}$ & 0.0 & \multirow{2}{*}{0.5} & \multirow{2}{*}{0.1} \\
        & R:& 0.125 & 0.0 & 0.0 & 0.0 & 0.1 & 0.75$\sqrt{\mu_0}$ & - $\sqrt{\mu_0}$ & 0.0 & & \\\hline
        \multirow{2}{*}{RP2} & L: & 1.08 & 1.2 & 0.01 & 0.5 & 0.95 & 2.0 & 3.6 & 2.0 & \multirow{2}{*}{0.4} & \multirow{2}{*}{0.2} \\
        & R:& 0.9891 & -0.0131 & 0.0269 & 0.010037 & 0.97159 & 2.0 & 4.0244 & 2.0026 & & \\\hline
        \multirow{2}{*}{RP3} & L: & 1.7 & 0.0 & 0.0 & 0.0 & 1.7 & 3.899398 & 3.544908 & 0.0 & \multirow{2}{*}{0.4} & \multirow{2}{*}{0.15} \\
        & R:& 0.2 & 0.0 & 0.0 & -1.496891 & 0.2 & 3.899398 & 2.785898 & 2.192064 & & \\\hline
        \multirow{2}{*}{RP4} & L: & 1.0 & 0.0 & 0.0 & 0.0 & 1.0 & 1.3$\sqrt{\mu_0}$ & $\sqrt{\mu_0}$ & 0.0 & \multirow{2}{*}{0.5} & \multirow{2}{*}{0.16} \\
        & R:& 0.4 & 0.0 & 0.0 & 0.0 & 0.4 & 1.3$\sqrt{\mu_0}$ & - $\sqrt{\mu_0}$ & 0.0 & & \\\hline
    \end{tabular}
    \caption{Initialization of Riemann Problems for the Euler equations. The initial conditions are given in terms of left (L) and right (R) states separated by an initial discontinuity at $x = x_0$ on the computational domain $[0,1]$. Further, the final time is denoted by $t_f$ and in all cases $\gamma = 5/3$ and $\mu_0 = 4\pi$. }
    \label{tab:RPMHD}
\end{table}

\subsubsection{Field Loop advection}

\begin{figure}[htpb]
    \centering
    \includegraphics[width=.45\textwidth]{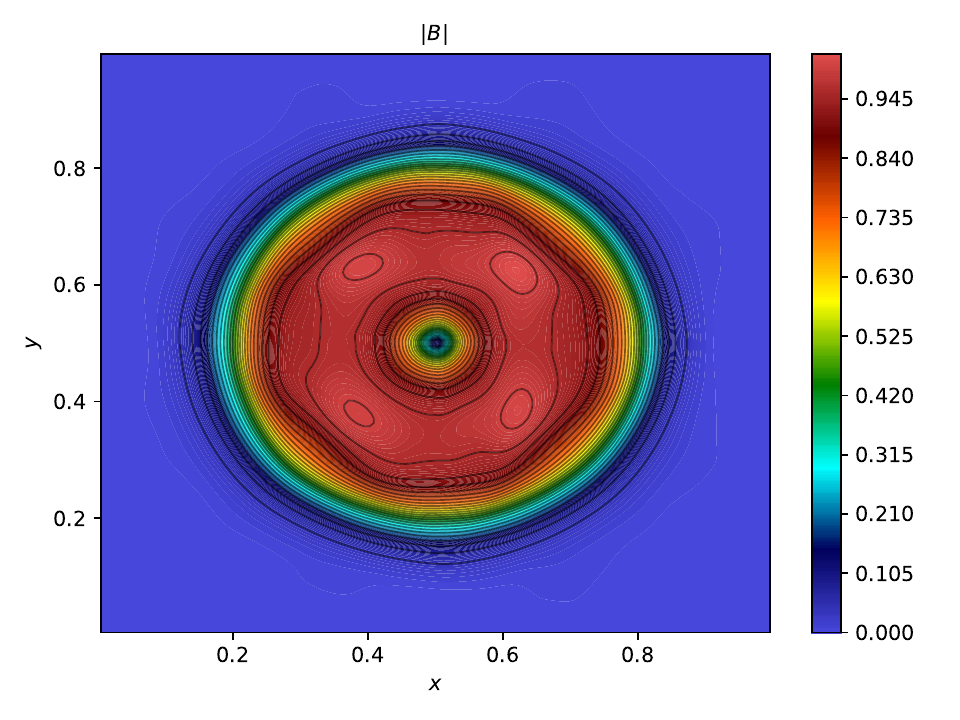}
    \includegraphics[width=.45\textwidth]{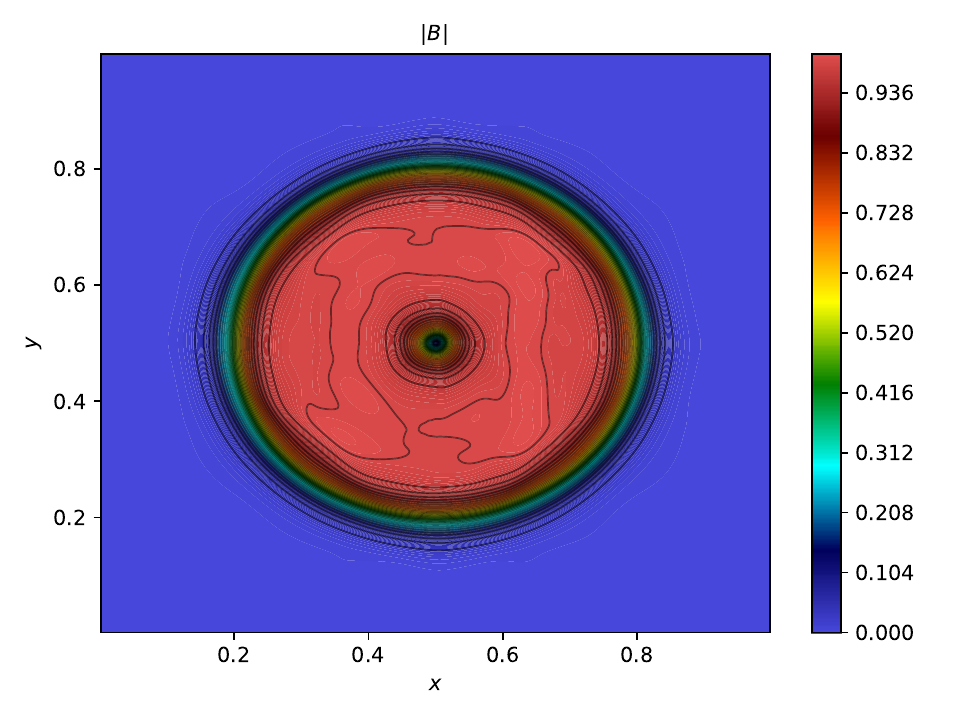}\\
    \includegraphics[width=.45\textwidth]{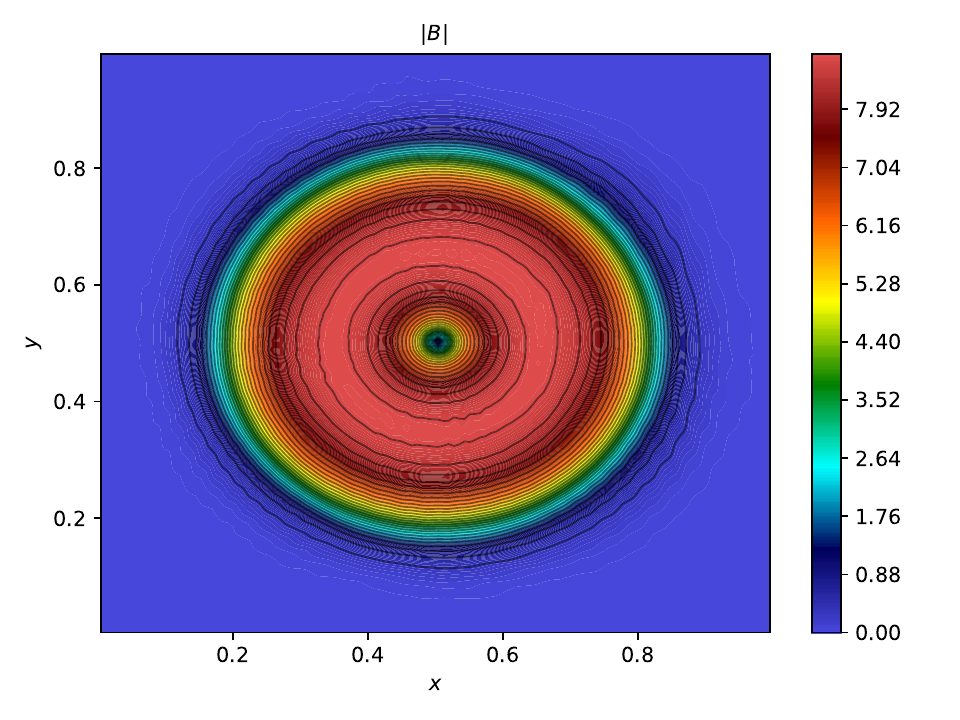}
    \includegraphics[width=.45\textwidth]{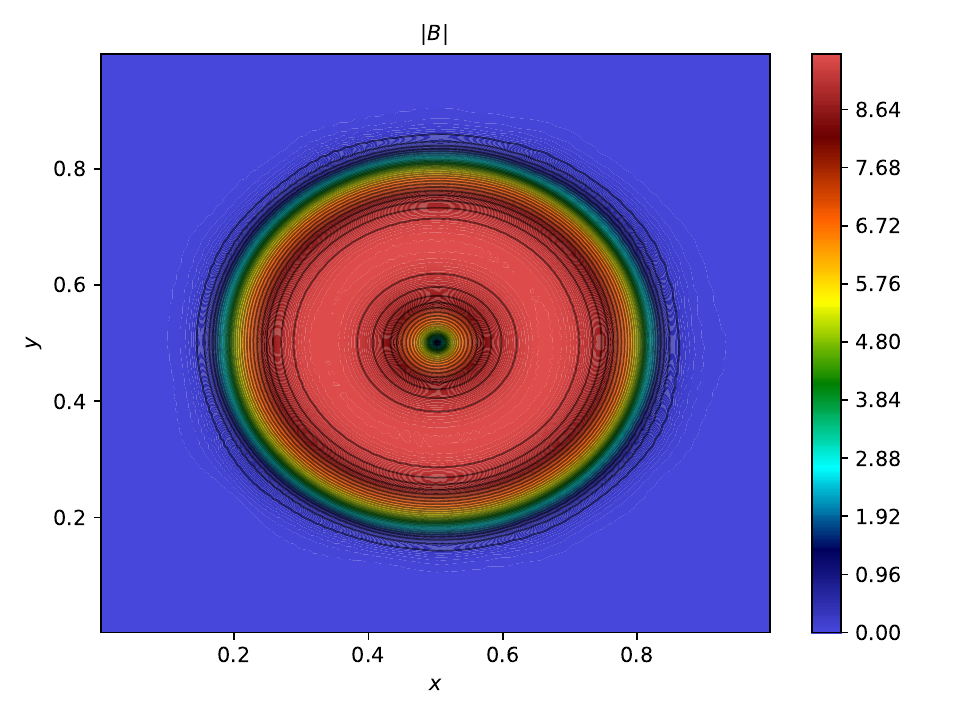}
    \caption{Field loop advection test at $t_f = 1$: Magnitude of the magnetic field with $A_0 = 1$ in the top row and $A_0 = 10$ in the bottom row computed on $N\times N$ cells displayed with 31 equidistant contour lines. Left collumn: $N=128$, Right column: $N=256$.  }
    \label{fig:FieldLoopAdvection}
\end{figure}

We solve the magnetic field loop advection problem originally proposed in \cite{Gardiner2005}, in the low Mach number version of \cite{Dumbser2019} with $\rho = 1$, an initial velocity field $(u_1,u_2) = (2,1)$ and the pressure $p = 10^5$. The computational domain is given by $\Omega=[-1;1]\times[-0.5;0.5]$ with periodic boundaries, and we use a mesh composed of $N \times N = 512 \times 256$ cells. The magnetic field is prescribed by the magnetic vector potential in $z$-direction given by
\begin{equation}
    A_z = \begin{cases}
        A_0 (R- r) & \text{if} \quad r \leq R,\\
        0          & \text{else},
    \end{cases}
\end{equation}
where $R = 0.3$ and $r = \sqrt{x_1^2 + x_2^2}$. This test case is very challenging because of a singular point in the magnetic field $\BB$ at the center of the computational domain.
In the setting of \cite{Dumbser2019} $A_0 = 10^{-3}$ which places the flow in a low acoustic Mach number regime $M_c \approx 6 \cdot 10^{-3}$ and high Alfv\'en Mach number regime $M_a \approx 7.9 \cdot 10^3$. To make the test more challenging, we perform four different simulations with $A_0 = 1, 10$, reaching a configuration with an Alfv\'en Mach number $M_a \approx 7.9 \cdot 10^{-1}$ corresponding to $A_0=10$. The final time is chosen to be $t_f = 1$, and the results are depicted in Figure \ref{fig:FieldLoopAdvection} for $|\BB|$.
The shape of the magnitude of the magnetic field remains qualitatively unchanged and shows only marginal influence on the diffusion despite the difference in the Alfv\'en Mach number, hence demonstrating that the new \IMEX scheme is stable independently of the magnetic scales. Indeed, the time step for all the simulations is the same and is only determined by the convective velocity $(u_1,u_2)=(2,1)$ which transports the loop diagonally through the domain.
In Figure \ref{fig:FieldLoopAdvection}, we show the results for two meshes made of $N\times N$ cells with $N = 128,256$ where the higher resolution yields a better capturing of the magnitude of the magnetic field.

\subsubsection{Orszag-Tang vortex}

\begin{figure}[htpb]
    \includegraphics[width=0.5\textwidth]{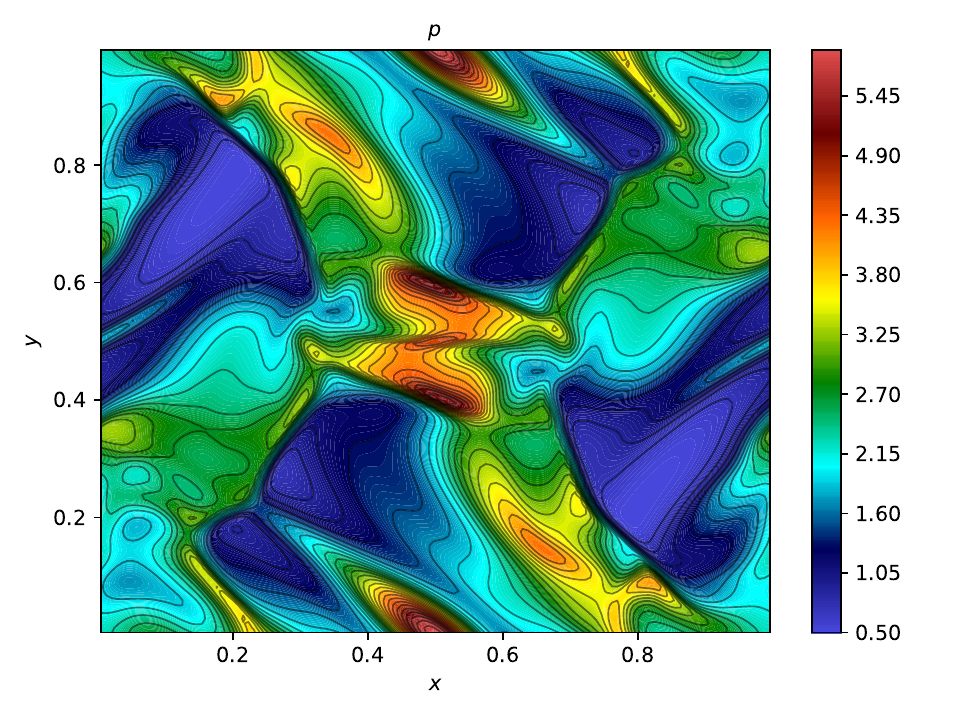}
    \includegraphics[width=0.5\textwidth]{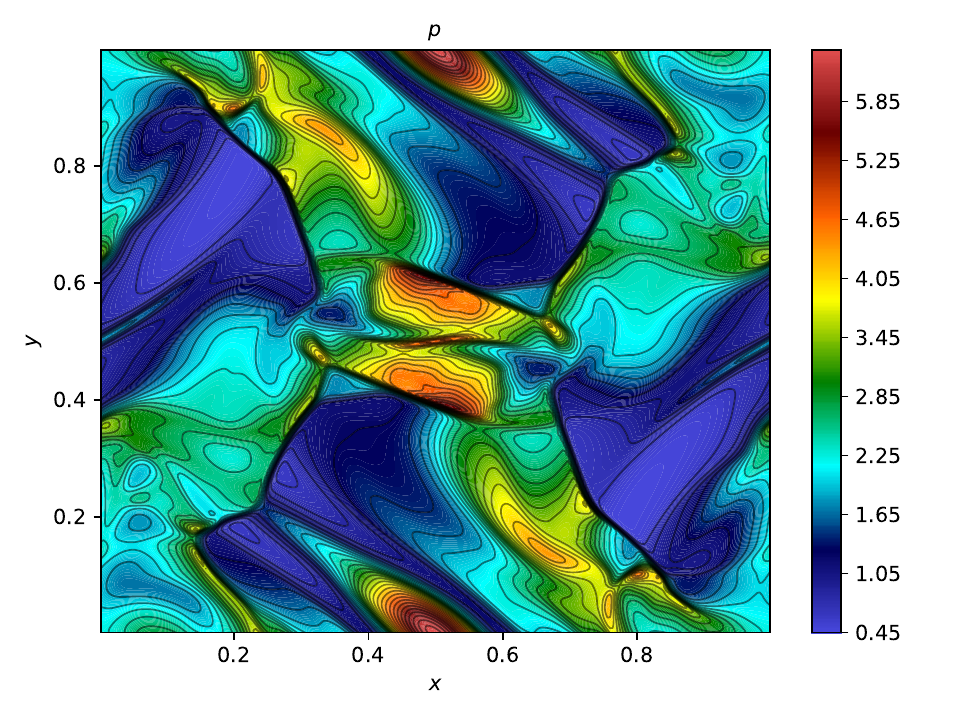}
    \caption{Orszag-Tang vortex: Pressure at time $t = 0.5$ computed on $N\times N$ cells displayed with 31 equidistant contour lines. Left: $N=128$, Right: $N=256$. }
    \label{fig:OrszagTang}
\end{figure}

We consider the well-known Orszag-Tang vortex for the two-dimensional ideal MHD equations \cite{OrszagTang1979,Dahlburg1989,Picone1991}.
The initial condition is smooth and reads
\begin{equation}
    (\rho, u_1, u_2, p, B_1, B_2) =
    \left(\frac{25}{9}, - \sin(2 \pi y), \sin(2 \pi x), \frac{5}{3}, - \sin (2 \pi y)\sqrt{4 \pi} , \sin (4 \pi y)\sqrt{4 \pi} \right).
\end{equation}
Though initially smooth, the dynamics develop shocks along the diagonal direction in combination with a vortex located at the center of the computational domain, which is defined by $\Omega=[0;1]\times [0;1]$ with periodic boundaries. We use a mesh made of $N\times N$ cells with $N=128, 256$. Figure \ref{fig:OrszagTang} depicts the pressure computed with the new second-order \IMEX scheme at times $t = 0.5$ for the two resolutions.
The features of the flows and the development of the shock are captured well and the results are qualitatively in very good agreement with those reported in the literature \cite{Dumbser2019,Balsara2015,Balsara2010,Boscheri2024}.

\section{Conclusions}
We have proposed a semi-implicit relaxed solver based on the Jin-Xin relaxation applied to stiff sub-fluxes for hyperbolic multi-scale systems of conservation laws.
The scheme was constructed in the cases where one or two stiff sub-systems are present in the model, i.e. where one or two scaling parameters determine the evolution of the flow.
In the considered models, it was the Mach number in Eulerian gas dynamics and the Mach and Alfv\'en number in the ideal MHD equations on which a two- and three-split scheme was applied respectively.
Due to the relaxation approach the resulting relaxation model was linear with respect to the stiff relaxation variables and the non-linear scale dependent sub-fluxes where shifted to algebraic source terms.
This structure was used to derive a prediction-correction method, where first a solution depending on the frozen characteristics of the relaxation model was obtain which was then corrected by projection on the equilibrium solution as the relaxation rate tended to zero.
To reduce computational overhead arising from the added relaxation parameters, the scheme was reformulated into independent wave-type equations which can be solved efficiently.
In case of the Euler equations due to the employed splitting the scheme was proven to be asymptotic preserving and contact preserving which was validated by numerical tests.
Within the framework of IMEX Runge-Kutta methods an extension to higher orders of accuracy was established.
In case of the three-split scheme applied to the ideal MHD equations, a second order extension could be achieved which was verified numerically.
However, this does not imply a general extension to higher order employing partitioned Runge-Kutta methods in time and could depend on the splitting.
Further investigations in the extension to more than three sub-systems and higher order extensions need to be done which is subject to future work.


\bibliographystyle{plain}
\bibliography{JinXin_multisplit}

\begin{thebibliography}{10}

\bibitem{AbbIolPup2017}
E.~Abbate, A.~Iollo, and G.~Puppo.
\newblock An all-speed relaxation scheme for gases and compressible materials.
\newblock {\em J. Comput. Phys.}, 351:1--24, 2017.

\bibitem{AvgBerIolRus2019}
S.~Avgerinos, F.~Bernard, A.~Iollo, and G.~Russo.
\newblock Linearly implicit all mach number shock capturing schemes for the
  euler equations.
\newblock {\em J. Comput. Phys.}, 393:278--312, 2019.

\bibitem{Balsara2004}
D.~Balsara.
\newblock {Second-Order Accurate Schemes for Magnetohydrodynamics with
  Divergence-Free Reconstruction}.
\newblock {\em Astrophys. J. Suppl. Ser.}, 151:149--184, 2004.

\bibitem{Balsara2010}
D.~S. Balsara.
\newblock {Multidimensional HLLE Riemann solver: application to Euler and
  magnetohydrodynamic flows}.
\newblock {\em J. Comput. Phys.}, 229(6):1970--1993, 2010.

\bibitem{Balsara2015}
D.~S. Balsara and M.~Dumbser.
\newblock {Divergence-free MHD on unstructured meshes using high order finite
  volume schemes based on multidimensional Riemann solvers}.
\newblock {\em J. Comput. Phys.}, 299:687--715, 2015.

\bibitem{Bispen2017}
G.~Bispen, M.~Luk{\'a}{\v{c}}ov{\'a}-Medvid'ov{\'a}, and L.~Yelash.
\newblock Asymptotic preserving {IMEX} finite volume schemes for low {M}ach
  number {E}uler equations with gravitation.
\newblock {\em J. Comput. Phys.}, 335:222--248, 2017.

\bibitem{Boscarino2016}
S.~Boscarino, F.~Filbet, and G.~Russo.
\newblock {High Order Semi-implicit Schemes for Time Dependent Partial
  Differential Equations}.
\newblock {\em J. Sci. Comput.}, 68:975--1001, 2016.

\bibitem{BoscarinoRussoScandurra2018}
S.~Boscarino, G.~Russo, and L.~Scandurra.
\newblock All {M}ach number second order semi-implicit scheme for the {E}uler
  equations of gas dynamics.
\newblock {\em J. Sci. Comput.}, 77(2):850--884, 2018.

\bibitem{Boscheri2020}
W.~Boscheri, G.~Dimarco, R.~Loub{\`{e}}re, M.~Tavelli, and M.H. Vignal.
\newblock {A second order all Mach number IMEX finite volume solver for the
  three dimensional Euler equations}.
\newblock {\em J. Comp. Phys.}, 415:109486, 2020.

\bibitem{Boscheri2021}
W.~Boscheri and L.~Pareschi.
\newblock {High order pressure-based semi-implicit IMEX schemes for the 3D
  Navier-Stokes equations at all Mach numbers}.
\newblock {\em J. Comput. Phys.}, 434:110206, 2021.

\bibitem{Boscheri2024}
W.~Boscheri and A.~Thomann.
\newblock {A structure-preserving semi-implicit IMEX Finite Volume Scheme for
  Ideal Magnetohydrodynamics at all Mach and Alfvén Numbers}.
\newblock {\em J. Sci. Comput.}, 100(3), July 2024.

\bibitem{CoulFraHelRatSon2019}
D.~Coulette, E.~Franck, P.~Helluy, A.~Ratnani, and E.~Sonnendr\"ucker.
\newblock Implicit time schemes for compressible fluid models based on
  relaxation methods.
\newblock {\em Comput. Fluids}, 188:70--85, 2019.

\bibitem{Cravero2018}
I.~Cravero, G.~Puppo, M.~Semplice, and G.~Visconti.
\newblock {CWENO: uniformly accurate reconstructions for balance laws}.
\newblock {\em Math. Comput.}, 87(312):1689--1719, 2018.

\bibitem{Dahlburg1989}
R.~B. Dahlburg and J.~M. Picone.
\newblock {Evolution of the Orszag--Tang vortex system in a compressible
  medium. I. Initial average subsonic flow}.
\newblock {\em Phys. Fluids B}, 1(11):2153--2171, 1989.

\bibitem{Brauer2016}
A.~de~Brauer, A.~Iollo, and T.~Milcent.
\newblock A {C}artesian scheme for compressible multimaterial models in 3d.
\newblock {\em J. Comput. Phys.}, 313:121--143, 2016.

\bibitem{Degond2011}
P.~Degond and M.~Tang.
\newblock {All speed scheme for the low {M}ach number limit of the isentropic
  {E}uler equations}.
\newblock {\em Commun. Comput. Phys.}, 10(1):1--31, 2011.

\bibitem{Dellacherie2010}
S.~Dellacherie.
\newblock Analysis of {G}odunov type schemes applied to the compressible
  {E}uler system at low {M}ach number.
\newblock {\em J. Comput. Phys.}, 229(4):978--1016, 2010.

\bibitem{Dumbser2016}
M.~Dumbser and D.~S. Balsara.
\newblock {A new efficient formulation of the HLLEM Riemann solver for general
  conservative and non-conservative hyperbolic systems}.
\newblock {\em J. Comput. Phys.}, 304:275--319, 2016.

\bibitem{Dumbser2019}
M.~Dumbser, D.S. Balsara, M.~Tavelli, and F.~Fambri.
\newblock A divergence-free semi-implicit finite volume scheme for ideal,
  viscous, and resistive magnetohydrodynamics.
\newblock {\em Int. J. Numer. Methods Fluids}, 89:16--42, 2019.

\bibitem{Dumbser2011}
M.~Dumbser and E.~F. Toro.
\newblock On universal osher-type schemes for general nonlinear hyperbolic
  conservation laws.
\newblock {\em Commun. Comput. Phys.}, 10(3):635--671, 2011.

\bibitem{Falle2002}
S.~A. E.~G. Falle.
\newblock {Rarefaction Shocks, Shock Errors, and Low Order of Accuracy in
  ZEUS}.
\newblock {\em ApJ}, 577(2):L123--L126, October 2002.

\bibitem{Falle2001}
S.~A. E.~G. Falle and S.~S. Komissarov.
\newblock On the inadmissibility of non-evolutionary shocks.
\newblock {\em J. Plasma Phys.}, 65(1):29--58, 2001.

\bibitem{Fambri2021}
F.~Fambri.
\newblock A novel structure preserving semi-implicit finite volume method for
  viscous and resistive magnetohydrodynamics.
\newblock {\em Int. J. Numer. Methods Fluids}, 93:3447--3489, 2021.

\bibitem{Feireisl2021}
E.~Feireisl, M.~Lukáčová-Medviďová, H.~Mizerová, and B.~She.
\newblock {\em Numerical Analysis of Compressible Fluid Flows}.
\newblock Springer International Publishing, 2021.

\bibitem{Gardiner2005}
T.A. Gardiner and J.M. Stone.
\newblock {An unsplit Godunov method for ideal MHD via constrained transport}.
\newblock {\em J. Comput. Phys.}, 205(2):509--539, 2005.

\bibitem{GuillardViozat1999}
H.~Guillard and C.~Viozat.
\newblock On the behaviour of upwind schemes in the low {M}ach number limit.
\newblock {\em Comput. Fluids}, 28(1):63--86, 1999.

\bibitem{Hairer1991}
E.~Hairer and G.~Wanner.
\newblock {\em Solving Ordinary Differential Equations II: Stiff and
  Differential-Algebraic Problems}.
\newblock Springer, New York, 1991.

\bibitem{JinXin1995}
S.~Jin and Z.~Xin.
\newblock The relaxation schemes for systems of conservation laws in arbitrary
  space dimensions.
\newblock {\em Commun. Pur. Appl. Math.}, 48(3):235--276, 1995.

\bibitem{KlainermanMajda1981}
S.~Klainerman and A.~Majda.
\newblock Singular limits of quasilinear hyperbolic systems with large
  parameters and the incompressible limit of compressible fluids.
\newblock {\em Commun. Pur. Appl. Math.}, 34(4):481--524, 1981.

\bibitem{Klein1995}
R.~Klein.
\newblock Semi-implicit extension of a {G}odunov-type scheme based on low
  {M}ach number asymptotics {I}: One-dimensional flow.
\newblock {\em J. Comput. Phys.}, 121(2):213 -- 237, 1995.

\bibitem{MichelDansac2022}
V.~Michel-Dansac and A.~Thomann.
\newblock {TVD-MOOD schemes based on implicit-explicit time integration}.
\newblock {\em Appl. Math. Comput.}, 433:127397, 2022.

\bibitem{Munz2000}
C-D. Munz, P.~Omnes, R.~Schneider, E.~Sonnendr{\"u}cker, and U.~Voss.
\newblock {Divergence correction techniques for Maxwell solvers based on a
  hyperbolic model}.
\newblock {\em J. Comput. Phys.}, 161(2):484--511, 2000.

\bibitem{Oeffner2025}
P.~\"Offner, L.~Petri, and D.~Torlo.
\newblock {Analysis for Implicit and Implicit-Explicit ADER and DeC Methods for
  Ordinary Differential Equations, Advection-Diffusion and Advection-Dispersion
  Equations}.
\newblock {\em Appl. Numer. Math.}, 2025.

\bibitem{OrszagTang1979}
S.~A. Orszag and C-M. Tang.
\newblock Small-scale structure of two-dimensional magnetohydrodynamic
  turbulence.
\newblock {\em J. Fluid Mech.}, 90(1):129--143, 1979.

\bibitem{Picone1991}
J.~M. Picone and R.~B. Dahlburg.
\newblock {Evolution of the Orszag--Tang vortex system in a compressible
  medium. II. Supersonic flow}.
\newblock {\em Phys. Fluids B}, 3(1):29--44, 1991.

\bibitem{Shu1998}
C.-W. Shu.
\newblock Essentially non-oscillatory and weighted essentially non-oscillatory
  schemes for hyperbolic conservation laws.
\newblock In {\em Advanced numerical approximation of nonlinear hyperbolic
  equations}, pages 325--432. Springer, 1998.

\bibitem{Sod1978}
G.~A. Sod.
\newblock A survey of several finite difference methods for systems of
  nonlinear hyperbolic conservation laws.
\newblock {\em J. Comput. Phys.}, 27(1):1--31, 1978.

\bibitem{Thomann2023}
A.~Thomann, A.~Iollo, and G.~Puppo.
\newblock Implicit relaxed all mach number schemes for gases and compressible
  materials.
\newblock {\em SIAM J. Sci. Comput.}, 45(5):A2632--A2656, 2023.

\bibitem{Toro2009}
E.~F. Toro.
\newblock {\em {Riemann Solvers and Numerical Methods for Fluid Dynamics }}.
\newblock Springer-Verlag: Berlin, 2009.

\bibitem{Toro2012}
E.~F. Toro and M.~E. V{\'a}zquez-Cend{\'o}n.
\newblock {Flux splitting schemes for the Euler equations}.
\newblock {\em Comput. \& Fluids}, 70:1--12, 2012.

\bibitem{torrilhon2003}
M.~Torrilhon.
\newblock {Non-uniform convergence of finite volume schemes for Riemann
  problems of ideal magnetohydrodynamics}.
\newblock {\em J. Comput. Phys.}, 192(1):73--94, 2003.

\end{thebibliography}
\end{document}